\documentclass[a4paper]{amsart}
\usepackage{fullpage}
\usepackage[english]{babel}
\usepackage[utf8]{inputenc}
\usepackage{amsmath,amsfonts,amsthm,amssymb}

\DeclareFontFamily{U}{mathx}{}
\DeclareFontShape{U}{mathx}{m}{n}{<-> mathx10}{}
\DeclareSymbolFont{mathx}{U}{mathx}{m}{n}
\DeclareMathAccent{\widehat}{0}{mathx}{"70}
\DeclareMathAccent{\widecheck}{0}{mathx}{"71}
\DeclareMathAccent{\widetilde}{0}{mathx}{"72}

\usepackage{mathtools}
\usepackage{graphicx}
\usepackage{subcaption}
\usepackage{type1cm}
\usepackage{eso-pic}
\usepackage{color}
\usepackage{pgfplotstable}
\usepackage{pgfplots}
\pgfplotsset{compat=1.5}
\usepackage{caption}
\makeatother
\usepackage{enumitem}
\usepackage{hyperref}
\usepackage{placeins}
\usepackage{todonotes}
\usepackage{comment}
\makeatletter
\AtBeginDocument{%
  \expandafter\renewcommand\expandafter\subsection\expandafter
    {\expandafter\@fb@secFB\subsection}%
  \newcommand\@fb@secFB{\FloatBarrier
    \gdef\@fb@afterHHook{\@fb@topbarrier \gdef\@fb@afterHHook{}}}%
  \g@addto@macro\@afterheading{\@fb@afterHHook}%
  \gdef\@fb@afterHHook{}%
}
\makeatother
\usepackage{tikz,tikz-cd}
\usepackage{float}
\usetikzlibrary{patterns}
\usetikzlibrary{calc}
\usepackage[capitalize,noabbrev]{cleveref}
\newtheorem{theorem}{Theorem}[section]
\newtheorem{lemma}[theorem]{Lemma}
\newtheorem{remark}[theorem]{Remark}
\newtheorem{problem}[theorem]{Problem}

\newtheorem{corollary}[theorem]{Corollary}
\newtheorem{definition}[theorem]{Definition}
\newtheorem{proposition}[theorem]{Proposition}

\newcommand{\calC}{\mathcal{C}}

\newcommand\dx{\,\mathrm{d}x}

\newcommand{\frakA}{\mathfrak{A}}

\newcommand{\calW}{\mathcal{W}}
\newcommand{\calT}{\mathcal{T}}
\newcommand{\calI}{\mathcal{I}}
\newcommand{\Nb}{\mathsf{N}}

\newcommand{\btheta}{{\boldsymbol{\operatorname{\theta}}}}

\DeclareMathOperator{\NN}{\boldsymbol{\mathcal{N}}(\btheta)}
\DeclareMathOperator{\NNv}{\boldsymbol{\mathcal{N}}}
\newcommand{\bW}{{\boldsymbol{\operatorname{W}}}}
\newcommand{\bb}{{\boldsymbol{\operatorname{b}}}}

\newcommand{\vt}{v_{\btheta}(\bx)}
\newcommand{\vl}{v_{\btheta}^\ell}

\DeclareMathOperator{\Div}{div}

\DeclareMathOperator{\Span}{span}

\DeclareMathOperator{\diag}{diag}
\DeclareMathOperator{\rank}{rank}
\DeclareMathOperator{\supp}{supp}

\DeclareMathOperator*{\argmin}{argmin}

\newcommand{\bbetap}{\boldsymbol{\beta}_\bp}
\newcommand{\bbetapl}{{\boldsymbol{\beta}_\bp^\ell}}

\newcommand{\Id}{{\operatorname{Id}}}

\newcommand{\eqnorm}[1]{{\left\vert\kern-0.25ex\left\vert\kern-0.25ex\left\vert #1
    \right\vert\kern-0.25ex\right\vert\kern-0.25ex\right\vert}}

\newcommand{\R}{{\mathbb R}}

\newcommand{\bx}{{\boldsymbol x}}

\newcommand{\bz}{{\boldsymbol z}}

\newcommand{\balpha}{{\boldsymbol \alpha}}

\newcommand{\bp}{{\boldsymbol p}}

\newcommand{\rom}[1]{\uppercase\expandafter{\romannumeral #1\relax}}

\usepackage[dvipsnames]{xcolor}

\floatstyle{ruled}
\newfloat{algorithm}{tbp}{loa}
\floatname{algorithm}{Algorithm}
\floatstyle{plain}

\title{Neural numerical homogenization based on Deep Ritz corrections}

\author[]{M.~Elasmi$^\dagger$, F.~Krumbiegel$^\dagger$, R.~Maier$^\dagger$}
\address{${}^{\dagger}$ Institute for Applied and Numerical Mathematics, Karlsruhe Institute of Technology, Englerstr.~2, 76131 Karlsruhe, Germany}
\email{\{mehdi.elasmi,felix.krumbiegel,roland.maier\}@kit.edu}

\date{\textbf{\today}}

\begin{document}

\begin{abstract}
Numerical homogenization methods aim at providing appropriate coarse-scale approximations of solutions to (elliptic) partial differential equations that involve highly oscillatory coefficients. The localized orthogonal decomposition (LOD) method is an effective way of dealing with such coefficients, especially if they are non-periodic and non-smooth. It modifies classical finite element basis functions by suitable fine-scale corrections. In this paper, we make use of the structure of the LOD method, but we propose to calculate the corrections based on a Deep Ritz approach involving a parametrization of the coefficients to tackle temporal variations or uncertainties. Numerical examples {for a parabolic model problem} are presented to assess the performance of the approach.\bigskip

% 35K20 Initial-boundary value problems for second-order parabolic equations
% 65M60 Finite element, Rayleigh-Ritz and Galerkin methods for initial value and initial-boundary value problems involving PDEs
% 68T07 Artificial neural networks and deep learning
\end{abstract}
\maketitle

%----------------------------------------------------------
% Intro
%----------------------------------------------------------

\section{Introduction}\label{sec:intro}
Composite materials play an important role in modern technologies, as they enable the combination of complex properties such as high strength, lightweight, thermal resistance, and electrical conductivity- characteristics that are often unattainable with plain homogeneous materials. A notable example for their application can be found in batteries. 
An interesting question in this context is the study of the thermal behavior of the cell during the charging/discharging protocol, especially during fast charging, where increased voltages and currents can lead to thermal breakdowns if not properly monitored and controlled. 
In the most simple manner, the temperature in a battery cell can be modeled by a parabolic equation (similar to~\cite{queisser}). More precisely, we seek the heat distribution $u$ that solves 
\begin{equation}\label{eq:model}
	\partial_t u(t,x)-\operatorname{div}(a (t,x) \nabla u(t,x)) = f(t,x), \qquad x \in \Omega, \quad t \in [0,T], 
\end{equation}
where $\Omega$ is an open bounded domain, $T$ the final time, $a$ a given  heat conductivity coefficient, and $f$ a heat source. Note that $a$ is oscillatory in space and may change over time. %
In particular, battery cells usually consist of a large amount of layers that have very small widths $\varepsilon$. Several engineering systems share strong similarities with batteries in that they are multiscale thermal problems governed by the heat equation, where fine-scale phenomena determine macroscopic behavior. %
A first example is nuclear fuel rods, where microscale porosity, cracking, and fission gas bubbles inside fuel pellets strongly affect the effective thermal conductivity and thus the macroscopic temperature field, see~\cite{Lucuta1996}. %
Furthermore, one could consider power electronic modules, in which microscale degradation of solder layers and thermal interface materials alters thermal resistance and leads to localized hot spots at the device level, for more information refer to~\cite{Pop2006}. %
In~\cite{Kamal1974}, a third example is presented describing curing of thermoset polymers or composite materials, where microscale reaction kinetics and evolving microstructure control heat generation and effective thermal properties during macroscopic heat transfer. %
In all these systems, as in batteries, microscale evolution feeds back into the heat equation through temperature-dependent source terms and effective material parameters. In this work we focus on the examples of a battery, however, we expect that from the general assumption on the material coefficient that other multiscale problems can be handled similarly. 

The multiscale structure makes the simulation using standard methods, such as the finite element method (FEM), inherently expensive. If $\varepsilon$ is not resolved by the mesh parameter, the solution behaves badly even if only macroscopic approximation properties are of interest. Therefore, very fine meshes are necessary to obtain reasonable results (this is also known as the pre-asymptotic behavior of the FEM). Thus, reducing the computational complexity with {tailor-made} multiscale methods is highly beneficial. Under additional assumptions on the coefficient~$a$ (such as periodicity and scale separation in~\cite{Ves22}), the use of multiscale methods based on analytical homogenization theory are an appropriate and efficient choice. 
The method applied by~\cite{Ves22} is the heterogeneous multiscale method (HMM), which has first been introduced in the elliptic setting in \cite{EE03,EE05} (see also \cite{AbdEEV12}). The methodology has since then been applied to various parabolic and parabolic-type problems, see, e.g.,~\cite{MinZ07,AbdV12,Abd16,AbdH17}. It may also be used to tackle highly oscillatory coefficients in space and time and stochastic partial differential equations (PDEs), see~\cite{EckV23} and~\cite{AbdP12}, respectively. The HMM calculates a coarse system matrix by approximating the homogenized coefficient from analytical homogenization theory. The system is then used to compute a coarse approximation of the (homogenized) solution to the original PDE. 
Another approach is the multiscale finite element method, see \cite{HouW97,AllB05}, which corrects coarse basis functions by solving operator-adapted problems in each element of the mesh. The corresponding analysis is based on analytical homogenization theory as well, which again requires structural assumptions.

In practice, the structure of materials in batteries is more involved, specifically the so-called active material is usually strongly heterogeneous as well as non-periodic. This promotes the use of more general approaches, so-called numerical homogenization techniques that provably work under minimal assumptions on the coefficient. A starting point for further developments has been the variational multiscale method introduced in~\cite{Hug95,HugFMQ98}, which decomposes the problem into an equation for a coarse-scale and a fine-scale function, respectively, and appropriately includes fine-scale information into the coarse-scale equation. Further strategies include rough polyharmonic splines~\cite{OwhZB14}, where basis functions are constructed based on a {constrained} minimization of the PDE operator in the $L^2$-norm, or gamblets~\cite{Owh17,OwhS19}, where the energy of basis functions is minimized. Gamblets are introduced based on a game theoretic approach to multiscale problems and have been applied to parabolic PDEs in~\cite{OwhZ17}. Another popular approach is the multiscale FE2 method introduced in \cite{Fey99}. 
Other approaches under minimal assumptions are generalized (multiscale) finite element methods~\cite{BabL11,EfeGH13,MaS22,MaSD22} that build basis functions based on locally solving spectral problems. 
For a more in-depth overview of multiscale methods, we refer to~\cite{AltHP21} and the textbooks~\cite{OwhS19,MalP20}. 

The multiscale setting including space- and time-dependent coefficients poses difficult questions since the above-mentioned approaches require an update of the spatial approximation spaces at a relatively high cost. Methods that work around this bottleneck are thus vital and examples dealing with this question employing space-time approaches include an approach by~\cite{OwhZ07}, where the basis functions of the {ansatz} and test space are subject to a global harmonic transform, or the high-dimensional sparse finite element method in~\cite{TanH19}, which requires a clear scale separation of the coefficient. The construction of optimal local subspaces based on the generalized finite element method for space-time coefficients is analyzed in \cite{SchS22,SchST23}. Another space-time multiscale method based on the generalized multiscale finite element method is considered in~\cite{ChuELY18}. 

The proposed method in this paper is based on the localized orthogonal decomposition (LOD) method introduced in~\cite{MalP14,HenP13} (see also~\cite{MalP20}), that constructs an appropriate orthogonal decomposition of an ideal approximation space and remaining fine-scale functions. The ideal approximation space can be well-approximated by correcting classical coarse-scale finite element functions on a chosen scale~$H > \varepsilon$ with the solutions to local auxiliary problems (so-called \emph{corrections}). 
The LOD method provably works under very general assumptions on the coefficient $a$, which makes it fitting to our problem. The approach has previously been applied to parabolic PDEs in~\cite{MalP18} and to coupled elliptic-parabolic problems in~\cite{MalP17,AltCMPP20}. 
The LOD method is particularly well-suited for time-dependent problems, where the modified approximation space is built in an offline phase. A substantial speedup is then achieved in the online phase when iterating through time as the resulting system matrices are of moderate size only. The offline phase where the LOD basis functions are set up is the bottleneck of the method. If the coefficient changes in time, substantial re-computations of the approximation space are necessary such that the method loses its efficiency in the online phase. In~\cite{LjuMM22} the time-dependence was tackled by posing the LOD method in a space-time setting. Here, we stick to the idea of updating the approximation spaces but we avoid expensive re-computations by employing an appropriate trained artificial neural network (ANN) that realizes the corrections and is parametrized by the coefficient. This allows for quick updates once the coefficient has changed and an adjustment is necessary. Since the constraint of the local correction problems does not depend on the coefficient, it can be enforced exactly by a fixed linear projection of the network output. The training of the networks then reduces to the unconstrained minimization of the energy of the correction problems, whose minimizers are exactly the LOD corrections, and requires no data apart from a set of coefficients.
Such a strategy is also useful if the specific heat conductivity is not known a priori and only an approximate expected solution may be computed. Note that the possibility to quickly obtain modified basis functions for a plethora of coefficients eventually makes up for the training phase and therefore provides a valuable enhancement of the LOD method for the setting of space-time oscillatory coefficients. 

The paper is structured as follows. In Section~\ref{sec:homogenization}, we establish the parabolic model problem considered in this paper and introduce the basics of the LOD method as the spatial discretization technique. In Section~\ref{sec:ann}, we explain how to compute the Deep Ritz-based corrections, which are eventually used to set up a modified LOD approach. Finally, we present numerical examples in Section~\ref{sec:examples}. 

%----------------------------------------------------------
% Model
%----------------------------------------------------------

\section{Numerical Homogenization}\label{sec:homogenization}
\subsection{Model problem}\label{subsec:heat_model}
Recall that we {aim at (approximately) solving} parabolic PDEs of the form~\eqref{eq:model}. More precisely, for a bounded Lipschitz domain~$\Omega\subset\R^d$ ($d=1,2,3$) with boundary $\Gamma\coloneqq\partial\Omega$ and {a final time}~$T>0$, we seek a function~$u$ such that
\begin{equation}\label{eq:heat_pde}
	\begin{aligned}
		\partial_t u -\operatorname{div}(a\nabla u) &= f && \quad\mathrm{in}~\Omega\times (0,T],\\
		u &= g && \quad\mathrm{on}~\Gamma\times (0,T],\\
		u &= u_0 && \quad \mathrm{in}~\Omega\times \{0\},
	\end{aligned}
\end{equation}
with an appropriate source term $f$, boundary function $g$, and initial condition $u_0$, where the coefficient~$a\in L^\infty(0,T; L^\infty(\Omega))$ might vary in space and time and satisfies $0<\alpha\leq a{(t,x)}\leq \beta <\infty$ for {almost every}~$x\in\Omega$ {and~$t \in [0,T]$}. We note that~$a$ could also be matrix-valued but for simplicity we focus on scalar coefficients only. Here, we have coefficients in mind that vary in space on a fine scale $0<\varepsilon\ll 1$ and possibly in time as well. 
Furthermore, the coefficient $a$ could inherit some form of uncertainty in the sense that on a subset~$\omega\subset\Omega$ the coefficient might not be known exactly except for upper and lower bounds. 
In the following, we only consider the case $g\equiv 0$ for simplicity. The variational formulation of~\eqref{eq:heat_pde} then seeks~$u \in L^2(0,T;H_0^1(\Omega))\cap H^1(0,T; H^{-1}(\Omega))$ that solves
\begin{equation}\label{eq:heat_var}
	\begin{aligned}
		\langle \partial_t u, v\rangle + \frakA(u,v) &= \int_\Omega f v \dx
	\end{aligned}
\end{equation}
for all $v \in H_0^1(\Omega)$, where 
\begin{equation*}
\langle w,v\rangle\coloneqq w(v)
\end{equation*}
for $w\in H^{-1}(\Omega)$ and $v\in H^1_0(\Omega)$ denotes the corresponding duality pairing, and
\begin{equation}\label{eq:a_inner_product}
	\frakA(w,v) := \int_\Omega a \nabla w \cdot \nabla v \dx
\end{equation}
for $w,v\in H_0^1(\Omega)$ denotes the $a$-induced inner product on $H_0^1(\Omega)$. Note that $\frakA$ {may depend} on time through the coefficient $a$.

\subsection{{Coarse discretization spaces}}
{The current setting brings about many difficulties. A first one is that classical spatial  discretization schemes as, e.g., the {FEM} may perform badly if the oscillations that are present in the coefficient are not resolved. 
More precisely, if the FEM is used for the parabolic problem~\eqref{eq:heat_var}, where the mesh size~$H>\varepsilon$ is larger than the finest oscillations of the coefficient, then the FEM approximation is expected to yield large errors, see also~\cite{BabO00}. %

{The goal of numerical homogenization is to solve the PDE~\eqref{eq:heat_var} on coarse meshes with an approximate solution that yields a good accuracy independent of oscillations on finer scales. This is done by correcting coarse piece-wise finite element functions such that they fit to the problem at hand. This construction results in coarse system matrices, which makes the solution process more efficient. This is particularly important for time-dependent setting, where linear systems need to be solved repeatedly. }%
In {the following subsections}, we introduce the basics of the approach and {lay the foundation for modifications} based on a Deep Ritz ansatz in the {next} section. }

Let $\calT_H$ be a quasi-uniform, quadrilateral mesh with mesh size $H >\varepsilon$ and let $V_H$ be the conforming subspace of $H_0^1(\Omega)$ defined by
\begin{multline}\label{eq:def_VH}
	V_H\coloneqq\{v\in C^0(\Omega)\mid\forall K\in\mathcal{T}_H\colon v\vert_K\text{ is a polynomial of partial degree}\leq 1\},
\end{multline}
which is a subspace of~$H_0^1(\Omega)$. We denote with $\{\Lambda_j\}_{j=1}^{N_H}$, where $N_H=\dim V_H$, the classical nodal basis of $V_H$. Further, we need an appropriate quasi-interpolation operator onto $V_H$, which is introduced in the following. 
\begin{definition}[Patch]\label{def:patch}
	Let $\ell\geq 1$ and $S \subseteq \Omega$. Define the \emph{patch (of order $\ell$) around $S$} by the recursion~$\Nb^{\ell+1}(S) = \Nb(\Nb^{\ell}(S))$ with $\Nb^1(S) := \Nb(S)$ and
	\begin{equation*}
		\Nb(S) := \bigcup \bigl\{\overline{K} \in \calT_h\,\vert\, \overline{S} \,\cap\, \overline{K}\neq \emptyset\bigl\}.
	\end{equation*} 
\end{definition}
\begin{definition}[Local projective quasi-interpolation operator]\label{def:qInterpolation}
	The operator $\calI_H\colon L^2(\Omega) \to V_H$ is called \emph{local projective quasi-interpolation  operator}, if it is a projection, i.e.,  $\calI_H\circ\calI_H=\calI_H$, and stable in the sense of 
	\begin{equation}\label{eq:stabInterpolation}
		H^{-1}\|v-\calI_H v\|_{L^2(K)} + \|\nabla \calI_H v\|_{L^2(K)} \leq C \|\nabla v\|_{L^2(\Nb(K))}
	\end{equation}
	for all $K \in \calT_H$ and all $v\in H_0^1(\Omega)$, where the constant $C$ is independent of $H$.
\end{definition}

{For readability, we omit a construction of a suitable quasi-interpolation operator as it is not relevant for the overall construction. We note, however, that a suitable candidate is presented in~\cite[Ex.~3.1]{MalP20}, which is also employed in our numerical examples. }

\subsection{{Localized orthogonal decomposition}}\label{subsec:lod}
{In this section, we define the previously mentioned corrections of coarse basis functions~$\Lambda_{j}$. %
The issue is that the} space $V_H$ is not able to capture fine oscillations of functions that lie in the so-called \emph{fine-scale space} $\calW := \ker \calI_H \vert_{H_0^1(\Omega)}$. Therefore, the main idea of the LOD method is to suitably {adapt} the basis functions~$\Lambda_j$ of $V_H$ by functions in the space $\calW$. We obtain the fine-scale information by the {following correction operators.}

\begin{definition}[Localized Corrections]\label{def:locCorrOp}
    {Let~$\Lambda_j$ be a nodal basis functions of $ V_H $ and $K\in\mathcal{T}_H$ be an element with~$\supp(\Lambda_j)\cap K\neq\emptyset$. Further, let~$\ell\geq 1$ be given. Then, the \emph{localized element-wise correction}~$\calC_K^\ell\Lambda_j$ of the basis function~$\Lambda_j$ is given by the following {constrained} minimization problem, which seeks the unique minimizer }
    \begin{multline}\label{eq:mincorr}
    	\calC_K^\ell\Lambda_j=\argmin_{v\in H_0^1(\Nb^\ell(K))}\frac{1}{2}\int_{\Nb^\ell(K)}a \nabla v \cdot \nabla v \dx - \int_K a\nabla \Lambda_j \cdot \nabla v \dx\\
        \mathrm{subject~to}\quad \calI_Hv=0.
    \end{multline}
    {The (corrected) \emph{multiscale basis functions}~$\widetilde{\Lambda}_j$ are then given as }
    \begin{equation}\label{eq:corrBF}
    	\widetilde \Lambda_j^\ell=\Lambda_j - \sum_{K\in\calT_H}\calC_K^\ell\Lambda_j = \Lambda_j - \sum_{K\in\calT_H\colon\operatorname{supp}(\Lambda_j)\cap K\neq\emptyset}\calC_K^\ell\Lambda_j.
    \end{equation}
    Finally, the \emph{localized multiscale space} is defined by %
    \begin{equation}\label{eq:lodspace}
        \widetilde{V}_H^\ell = \Span\{\widetilde{\Lambda}_j^\ell\}_j^{N_H}. 
    \end{equation}
\end{definition}
We note here that due to the constraint in~\eqref{eq:mincorr} the corrections~$\calC_K^\ell\Lambda_j\in\mathcal{W}$ lie in the fine-scale space and contain the problem-adapted corrections that are not captured by coarse functions. Furthermore, there are many different presentations for the computation of the corrections, see also~\cite[Ch.~3.5]{MalP20}. %

{In practice, the minimization problem~\eqref{eq:mincorr} is solved using a fine finite element space on a small patch with mesh size~$h<\varepsilon$. %
We have that all corrector problems for each basis function~$\Lambda_j$ are independent of each other and can thus be computed embarrassingly parallel, which allows for a distribution of the computational resources across many cores. Further, for time-independent coefficients, the construction of the multiscale basis functions can be done a priori once in an offline phase and the computed basis functions can be used for every time step. This greatly decreases the online time as much smaller systems compared to a classical FEM need to be solved at each time. }

{However, in this work we specifically consider coefficients that change in time. Then the corrector problems have to be solved in (almost) every time step, which greatly increases the computational effort as we cannot make use of pre-computations, as the space~\eqref{eq:lodspace} becomes time-dependent. %
Possible solutions to overcome re-computations include using space-time ansatz functions as in~\cite{LjuMM22}. A disadvantage there is that we would need to compute all basis functions in the space-time mesh in advance, which introduces a large number of corrector problems and required storage. %
Another option is to use an error estimator to only partially update the spaces similar as in~\cite{HelM19,HelKM20,MaiV22}. However, if the coefficient is rapidly changing in the whole domain, all basis functions would need to be updated. %
{This motivates employing ANNs in the following section that act as a solution map from a given (local) coefficient to a corrected basis function, since ANNs (once trained) can be used online as the forward pass through a neural network is generally cheap. }%
As the corrections~\eqref{eq:mincorr} are defined based on an energy minimization problem, a Deep Ritz framework as introduced in~\cite{EY18} appears to be a natural choice. }

\begin{figure}
  \centering
	\begin{tikzpicture}[scale=0.37]
		% corner
		\fill[black!40!white] (0,0) rectangle (3,3);
		\fill[black!80!white] (0,0) rectangle (1,1);
		\node[white] at (0.5, 0.5) {$\pmb{K_2}$};

		% edge corner
		\fill[black!40!white] (11,10) rectangle (14,14);
		\fill[black!80!white] (13,12) rectangle (14,13);
		\node[white] at (13.5, 12.5) {$\pmb{K_5}$};

		\fill[black!40!white] (10,0) rectangle (14,4);
		\fill[black!80!white] (12,1) rectangle (13,2);
		\node[white] at (12.5, 1.5) {$\pmb{K_6}$};

		% edge center
		\fill[black!40!white] (0,6) rectangle (3,11);
		\fill[black!80!white] (0,8) rectangle (1,9);
		\node[white] at (0.5, 8.5) {$\pmb{K_1}$};

		\fill[black!40!white] (5,10) rectangle (10,14);
		\fill[black!80!white] (7,12) rectangle (8,13);
		\node[white] at (7.5, 12.5) {$\pmb{K_3}$};
		
		% center
		\fill[black!40!white] (4,3) rectangle (9,8);
		\fill[black!80!white] (6,5) rectangle (7,6);
		\node[white] at (6.5, 5.5) {$\pmb{K_4}$};
		
		% grid
		\draw[step=1.0,black,thin] (0,0) grid (14,14);
	
	\end{tikzpicture}
	\caption{Examples of {possible} patch configurations}\label{fig:patch}
\end{figure}

\subsection{{Fully discrete method}} 
{Finally, we present a fully discrete version of the LOD method applied to the parabolic model problem. Let $\widetilde u_H^m\in \widetilde V_H^\ell$ be the LOD solution and ${u_h^m\in V_h}$ be the fine finite element solution, respectively, at time point $m$ if~\eqref{eq:heat_var} is discretized in time with a backward Euler scheme. That is,}
{\begin{subequations}\label{eq:fd_heat}
\begin{align}
	% \begin{aligned}
		\int_{\Omega}\tfrac{\widetilde u_H^m - \widetilde u_H^{m-1}}{\tau} \, \widetilde v_H\dx + \frakA(\widetilde u_H^m, \widetilde v_H)& =\int_{\Omega}f \widetilde v_H \dx,\label{eq:fd_heat_lod}\\
		\int_{\Omega}\tfrac{u_h^m - u_h^{m-1}}{\tau} \, v_h\dx + \frakA(u_h^m, v_h) & = \int_{\Omega} f v_h\dx, \label{eq:fd_heat_fem}
	% \end{aligned}
\end{align}
\end{subequations}}
for all $\widetilde v_H\in\widetilde V_H^\ell$ and for all ${v_h\in V_h}$. %
{Note that, independently of the fine scale $\varepsilon$, the LOD solution $\widetilde{u}_H^m$ of \eqref{eq:fd_heat_lod} is an appropriate coarse-scale approximation of the solution of \eqref{eq:fd_heat_fem}, see, e.g., \cite[Ch.~9]{MalP20}.}  
%

%----------------------------------------------------------
% ANN
%----------------------------------------------------------
\section{Approximation of the Corrections with Artifical Neural Networks}\label{sec:ann}

{In this section, we describe how the corrections~\eqref{eq:mincorr} are approximated by ANNs. The coefficient is described by a parameter vector~$\bp$, and for each coarse element~$K$, the corrections of the basis functions are represented by an ANN with the input $(x,\bp)$. The ANN is trained by minimizing the energy in~\eqref{eq:mincorr} for a set of parameter samples, i.e., we follow a Deep Ritz approach. The training only requires the coefficients for these samples, but no precomputed corrections. Once trained, the networks provide the corrections, and thus the LOD basis, for any parameter by a forward pass, which allows for fast updates of the multiscale space in time-dependent or uncertain settings. We refer to the resulting method as Deep Ritz LOD (DR-LOD).}

\subsection{{Artificial neural networks}}
{We approximate the corrections by feed-forward ANNs, also known as multilayer perceptrons, which consist of a finite number of fully connected layers~\cite{goodfellow2016deep}. Each layer applies an affine transformation followed by a component-wise non-linear activation function~\cite{pinkus1999approximation}. An ANN of depth~$D>0$ with layer widths $\{N_k\}_{k=0}^{D}$ is characterized by its trainable parameters
\begin{equation*}
	\btheta := \{(\bW_k,\bb_k)\}_{k=1}^{D}, \qquad \bW_k \in \R^{N_k\times N_{k-1}},\quad \bb_k \in \R^{N_k},
\end{equation*}
i.e., the weight matrices and bias vectors of the layers, with $$N_\Theta := \sum_{k=1}^{D} N_{k}(N_{k-1}+1)$$ parameters in total.\\}
Next, we introduce the realization function for the type of ANNs considered in this work. 

\begin{definition}[Realization and network function]\label{def:RealFunction}
Let $\btheta:= \{(\bW_k,\bb_k)\}_{k=1}^{D}$, which defines an ANN of depth~$D>0$ and widths $\{N_k\}_{k=0}^{D}$. Moreover, let $\bx \in \R^{N_0}$ denote an input of the ANN. The realization function~$\NNv:\R^{N_\Theta} \to V_\btheta$ yields a network function~$v_\btheta:= \NN$ mapping the input $\bx$ to an output~$\vt \in \R^{N_D}$ such that 
\begin{equation*}
  \vt := \NN(\bx) = \mathcal{N}^{(D)} \circ \mathcal{N}^{(D-1)} \circ \dots \circ \mathcal{N}^{(1)}\circ \mathcal{N}^{(0)}(\bx),
\end{equation*}
where the realization functions of the different layers, denoted by $\mathcal{N}^{(\cdot)}$, are defined as 
\begin{align*}
	x_0&=\mathcal{N}^{(0)}(\bx) = \bx, \\
	x_k&=\mathcal{N}^{(k)}(\bx_{k-1}) = \sigma (\bW_k \bx_{k-1} + \bb_k), \qquad k = 1,\dots,D-1, \\
	\vt&=\mathcal{N}^{(D)}(\bx_{D-1}) = \bW_D \bx_{D-1} + \bb_D.
\end{align*}   
Here, $\sigma: \R \to I \subset \R$ is a possibly non-linear activation function and $V_\btheta \subset H^1(\Omega)$ denotes the space spanned by all possible network realizations. 
\end{definition}
{We use $\sigma := \tanh\colon \R \to (-1,1)$ as activation function, which is non-linear, bounded, and smooth. The trainable parameters are initialized by a Xavier/Glorot normal initialization~\cite{glorot2010understanding} to prevent exploding or vanishing gradients, i.e., the entries of $\bW_k$ and $\bb_k$ are drawn from normal distributions with zero mean and variance $\frac{2}{N_k + N_{k-1}}$ and $\frac{1}{N_k}$, respectively, $k = 1,\dots,D$.}
{
\begin{remark}[Choice of the architecture]
{In this study, we use simple fully connected networks. The main difficulty for such networks is the fine oscillation of the coefficient: by the frequency principle (F-principle)~\cite{xu2020frequency}, see also~\cite{RahBADLHBC19}, neural networks learn low frequencies first and high frequencies only slowly, which is why a single network on the whole domain struggles with rapidly varying coefficients. The LOD structure reduces this difficulty, since each network only has to represent a correction on a patch, which is mapped onto $[-1,1]^d$, see Section~\ref{subsec:trainingDetails}. Relative to the input domain of the network, the oscillations are therefore coarser by the factor $1/((2\ell+1)H)$. More expressive architectures such as DeepONets or Fourier neural operators could be combined with our approach in the same way.}
\end{remark}
}

\subsection{{Loss functions and discretization}}
{As mentioned in Section~\ref{subsec:lod}, we follow a Deep Ritz approach, i.e., the loss function is given by the energy of the minimization problem~\eqref{eq:mincorr}, in contrast to physics-informed neural networks (PINNs)~\cite{raissi2019physics}, which minimize the residual of the governing equations. Apart from a set of parametrized coefficients, no data is needed for the training. To replace $v \in H_0^1(\Nb^\ell(K))$ in~\eqref{eq:mincorr} by a network function, we first transform the output of the network such that it satisfies the homogeneous Dirichlet boundary condition on the patch.}
\begin{definition}[Localization of the network function]\label{def:locNF}
	Let the realization function $\NNv$ define an ANN with fixed architecture according to Definition~\ref{def:RealFunction}, and let $v_\btheta:=\NN \in V_\btheta \subset H^1(\Omega)$ denote the resulting network function. For $\Omega^\ell \subset \Omega \subseteq \R^d$, $d\geq1$, and $x \in \Omega$, a localized network function $\vl \in V_\btheta(\Omega^\ell)$ is understood as   
	\begin{equation*}
		\vl(x) = \begin{cases}
			g(v_\btheta(x)) & \text{ for } x \in \Omega^\ell, \\
			0 & \text{ otherwise}, 
		\end{cases}
	\end{equation*}
	where $g: \R^{N_D} \to \R^{N_D}$ is a smooth function that vanishes on $\partial \Omega^\ell$ and is non-zero in $\Omega^\ell$. {For $\Omega^\ell:= (x^a_1,x^b_1) \times \dots \times (x^a_d,x^b_d)$ and $x := (x_1,\dots,x_d)^\top \in \R^d$, the considered function $g$ is chosen as
	\[g(v_\btheta(x)) = \prod_{i=1}^d \frac{2}{\pi}\arctan\Big(\frac{1}{\kappa}\sin\Big(\pi\,\frac{x_i-x_i^a}{x_i^b-x_i^a}\Big)\Big)\, v_\btheta(x),  \]
    where $\kappa > 0$ is a small parameter. The function $g$ is a smoothed indicator function of the patch, i.e., it vanishes on $\partial\Omega^\ell$ and is close to one in the interior of $\Omega^\ell$.}
    Note that in practice only~$x \in \overline{\Omega^\ell}$ are taken into account during training. That is, we interpret $V_\btheta(\Omega^\ell)$ as a subspace of~$H^1_0(\Omega^\ell)$.
\end{definition}%

\begin{remark}
    {The concatenation of a neural network with a {function} as in Definition~\ref{def:locNF} forces the output of the network to fulfill the boundary condition on the patch exactly. This method is commonly used, see also~\cite[Sec.~5.1.1]{SukS22} for more details. }
\end{remark}

Without loss of generality, let us assume that the considered inputs of the ANN have the form~$\bx:=(x,\bp)$, where $x \in \Omega$ denotes the spatial coordinates and $\bp \in \R^{N_p}$, $N_p \geq 0$, refers to a possible additional parametrization vector for the coefficient function $a: \Omega \times \R^{N_p} \to \R$, cf. Section~\ref{subsec:lod}.
\begin{remark}\label{rem:inputsForm}
	Considering inputs of the form $\bx := (x,\bp)$ allows us to clearly separate the spatial component from the parametrization, which constitutes the variable of interest. 
\end{remark}

As a next step, we define the {loss function associated with the energy in~\eqref{eq:mincorr}.}
Afterwards, we introduce the fully discrete setting (including quadrature) that is used for computational purposes. 

\begin{definition}[{Loss function}]\label{def:lossFunction}
{Let $\ell \geq 1$ and $K \in \calT_H$ be given, and let $J(K):=\{j\,:\,K\subseteq\supp(\Lambda_j)\}$ denote the indices of the $2^d$ nodal basis functions associated with the vertices of~$K$.\footnote{For vertices on~$\partial\Omega$, $\Lambda_j$ denotes the nodal basis function of the space~\eqref{eq:def_VH} without boundary conditions. The corresponding corrections are not needed for homogeneous Dirichlet boundary conditions.} Further, for $j\in J(K)$, $v_{\btheta,j}^{\ell,K} \in V_\btheta(\Nb^\ell(K))$ denotes the component associated with~$\Lambda_j$ of a localized network function with the output components $(v^{\ell,K}_{\btheta,j})_{j\in J(K)}$ as given in Definition~\ref{def:locNF} with $\Omega^\ell = \Nb^\ell(K)$.}
Given the parametrization vector $\bp$, we define the energy loss function~$\mathcal{L}_\mathrm{energy}\colon V_\btheta(\Nb^\ell({K})) \to \R$ by
\begin{equation}\label{eq:energyLoss}
\begin{aligned}
	\mathcal{L}_\mathrm{energy}(v_{\btheta,j}^{\ell,{K}}) &= \frac{1}{2}\int_{\Nb^\ell({K})} a(x,\bp) \nabla v_{\btheta,j}^{\ell,{K}}(x,\bp) \cdot \nabla v_{\btheta,j}^{\ell,{K}}(x,\bp) \dx\\
	&\quad -\int_{{K}} a(x,\bp) \nabla \Lambda_j(x) \cdot \nabla v_{\btheta,j}^{\ell,{K}}(x,\bp)\dx{.}
\end{aligned}
\end{equation}
\end{definition}

\begin{remark}
    {The loss function in Definition~\ref{def:lossFunction} is a direct implication of the definition of the corrector problem as a minimization of the energy inner product (subject to a constraint). An advantage of this choice is that the implementation can be done straight-forward, avoiding the computation of the specifically constructed saddle-point problems, see, e.g.,~\cite{MalP20}. Furthermore, the quadrature of choice introduces a comparably small error of the fine-scale mesh that is negligible compared to the coarse mesh size where the solution is computed.}
\end{remark}

The computation of the integrals in~\eqref{eq:energyLoss} can be performed using numerical quadrature rules. In this case, the quadrature points would constitute the spatial component of the input. 
As stated in Remark~\ref{rem:inputsForm}, we are rather interested in the prediction of correction functions for arbitrary coefficients~$a$ in some given spectrum, which is determined by a set of distinct parametrization vectors collected in~${\bp}:= \{\bp_k\}_{k=1}^{N_\mathrm{s}}$, with~$N_\mathrm{s} \geq 1$. We consider a simple quadrature rule based on a classical space discretization of~\eqref{eq:energyLoss} using a finite element ansatz. More precisely, we consider a finite element space $V_h$ defined analogously to~\eqref{eq:def_VH} with mesh size $h<\varepsilon$ that resolves the fine scale. Further, we consider the underlying mesh $\calT_h$ to be a refinement of $\calT_H$, such that $V_H\subset V_h$. We denote with $\{\lambda_j\}_{j=1}^{N_h}$ the nodal basis of~$V_h$, where~$N_h = \dim V_h$. 
For convenience and as a preliminary step, we first consider global correction problems on $\Omega $, i.e., with a patch big enough to contain the whole domain ($\ell = \infty$). Based on the fine mesh $\calT_h$, the network function is approximated for each sample $\bp_k \in \R^{N_p}$ by
\begin{equation}\label{eq:nodalInterp}
	v_{\btheta,j}^{\infty,{K}}(x,\bp_k) = \sum_{{i}=1}^{N_h} {\beta}_{{i}}(z_{{i}},\bp_k) \lambda_{{i}}(x) , \quad x \in \Omega ,
\end{equation} 
where $\bz:=\{z_{{i}}\}_{{i}=1}^{N_h}$, $z_{{i}} \in \R^d${,} contains the nodes of the fine mesh $\calT_h$.
Hence, the output of the neural network reduces to $\boldsymbol{{\beta}}(\bz,\bp_k):= ({\beta}_1(z_1,\bp_k),\dots, \break{\beta}_{N_h}(z_{N_h},\bp_k))$. For notational simplicity, we omit the specification of the spatial component $\bz$ in the argument of ${\boldsymbol{\beta}}$ and write $\boldsymbol{{\beta}}_{\bp_k}$ instead of~$\boldsymbol{{\beta}}(\bz,\bp_k)$. Altogether, the input dataset for the ANN is defined by $\widecheck{\bx}:= \bz \times {\bp} \in \R^{N_0 \times N_h N_\mathrm{s}}$ with input dimension $N_0 = d +N_p$ and~$N_h N_\mathrm{s}$ training samples. The resulting output is denoted by~${\boldsymbol{\beta}}(\widecheck{\bx}):=\bbetap:= (\boldsymbol{{\beta}}_{\bp_1},\dots,\boldsymbol{{\beta}}_{\bp_{N_\mathrm{s}}})^\top \in \R^{N_h N_\mathrm{s}}$. 

To approximate~\eqref{eq:energyLoss}, let first for a given $j \in \{1,\ldots, N_H\}$ the vector $\boldsymbol{\gamma}:= (\gamma_1,\ldots,\gamma_{N_h})^\top \in \R^{N_h}$ be such that
\begin{equation*}
	\Lambda_j(x) = \sum_{{i}=1}^{N_h} \gamma_i \lambda_i(x)
\end{equation*}
represents the coarse basis function $\Lambda_j \in V_H$ via the fine basis functions of $V_h$.
Next, we define $\boldsymbol{\Lambda}_j:= \mathbf{1} \otimes \boldsymbol{\gamma}$ with $\textbf{1} = (1,\ldots,1)^\top \in \R^{N_\mathrm{s}}$ as the $N_\mathrm{s}$-times repeated vector $\boldsymbol{\gamma}$ (using the Kronecker product). Further, $\boldsymbol{S}^h_\bp:= \diag(S_h(\bp_1),\ldots,S_h(\bp_{N_\mathrm{s}}))$ denotes the block-diagonal matrix with $N_\mathrm{s}$ blocks, whose elements~$S_h(\bp_k) \in \R^{N_h \times N_h}$ are the stiffness matrices with respect to~the fine mesh $\calT_h$ and the parameter sample~$\bp_k$. We then obtain the following approximation of~\eqref{eq:energyLoss} 
\begin{equation}\label{eq:approxLossEnergy}
	 \widehat{\mathcal{L}}_\mathrm{energy}(\bbetap) = \frac{1}{2}\bbetap^\top \boldsymbol{S}_\bp^h \bbetap - \boldsymbol{\Lambda}_j^\top{\boldsymbol{S}_{\bp,K}^h} \bbetap{,}
\end{equation}
{where $\boldsymbol{S}_{\bp,K}^h$ is defined as $\boldsymbol{S}_\bp^h$ but with the stiffness matrices assembled on the element~$K$ only.}
Note that using the approximated loss function~\eqref{eq:approxLossEnergy} in the training procedure is a significant simplification compared to the original problem by avoiding the online computation of gradients and reducing the number of sensors, i.e., points where the solution should be learned, to the dimension of the fine mesh~$V_h$ for each $\bp_k$.

{Further, we denote by $I_H\colon \R^{N_h} \to \R^{N_H}$ the matrix version of the quasi-interpolation operator~$\calI_H$ for functions in~$V_h$, and by $\boldsymbol{I}_H = \diag(I_H,\ldots,I_H)$ the block-diagonal matrix with $N_\mathrm{s}$ blocks containing~$I_H$.} 

Similarly to Definition~\ref{def:locNF}, we address the localization procedure for the discrete setting and the adopted notations. Given a localization parameter~$\ell \geq 1$, we denote the localized coefficients corresponding to $v_{\btheta,j}^{\ell,{K}}$ by $\bbetapl$. {It is the restriction of~$\bbetap$ to the nodes} 
\begin{equation*} 
	\bz^\ell:=\{z_k\}_{k=1}^{N_h^\ell}
\end{equation*} 
containing all $ z_k \in \Nb^\ell({K})$, where $N_h^\ell = \dim V_h(\Nb^\ell({K})) \leq N_h$. Note that we implicitly assume for the moment that the first~$N_h^\ell$ nodes of the whole mesh lie in the patch $\Nb^\ell({K})$ to avoid a re-numbering. The coefficients associated to $z_k \notin \Nb^\ell({K})$ are fixed to zero and thus excluded from the training process. Correspondingly, we refer to the localized counterparts of $\boldsymbol{S}_\bp^h$, {$\boldsymbol{S}_{\bp,K}^h$,} $\boldsymbol{\Lambda}_j$, and $\boldsymbol{I}_H$ restricted to $\Nb^\ell({K})$ by $\boldsymbol{S}_\bp^{h,\ell}$, {$\boldsymbol{S}_{\bp,K}^{h,\ell}$,} $\boldsymbol{\Lambda}_j^\ell$, and $\boldsymbol{I}^\ell_H$, respectively.

\subsection{{Exact interpolation condition}}\label{subsec:exact}
{A straightforward way to incorporate the interpolation condition is to add a weighted penalty term $\|\boldsymbol{I}_H\bbetap\|_2^2$ to~\eqref{eq:approxLossEnergy}, where the weights may be determined in a self-adaptive way, following the idea of self-adaptive PINNs~\cite{mcclenny2023self}. However, this leads to a saddle-point problem, and the balance of the two loss terms is delicate during training, cf.~\cite{wang2022and}. On the other hand, the interpolation condition is linear and, in contrast to the energy, does not depend on the coefficient~$a$ or the parametrization~$\bp$. In the following, we exploit this and enforce the condition exactly by a fixed linear map applied to the network output, such that no penalty term and no weights are needed. We proceed in three steps: we characterize the admissible coefficient vectors of the localized corrections, construct a projection onto them, and show that minimizing the energy of projected vectors yields exactly the discrete LOD corrections. Throughout this subsection, we fix an element~$K$, a basis function~$\Lambda_j$ with $K\subseteq\supp(\Lambda_j)$, and a sample~$\bp_k$.}

\paragraph{{Admissible set}}
{Let $z_1,\dots,z_{n_0}$ denote the nodes in the interior of~$\Nb^\ell(K)$. We identify a vector $\boldsymbol{\beta}=(\beta_i)_{i=1}^{n_0}\in\R^{n_0}$, whose components are the values at these nodes, with the function $v_{\boldsymbol{\beta}}:=\sum_{i=1}^{n_0}\beta_i\lambda_{z_i}\in V_h\cap H^1_0(\Nb^\ell(K))$, where $\lambda_{z_i}$ denotes the nodal basis function of~$V_h$ associated with~$z_i$. This defines a bijection from~$\R^{n_0}$ onto $V_h\cap H^1_0(\Nb^\ell(K))$. In particular, the localized network output, which vanishes on~$\partial\Nb^\ell(K)$ by Definition~\ref{def:locNF}, is determined by its values at the interior nodes. Let $I_H^\ell$ denote the block of~$\boldsymbol{I}_H^\ell$ associated with one sample, and let $\mathbf{B}\in\R^{N_H^\ell\times n_0}$ consist of the columns of~$I_H^\ell$ that belong to the interior nodes. Since $v_{\boldsymbol{\beta}}$ vanishes at the remaining nodes of the patch, we have $\calI_H v_{\boldsymbol{\beta}}=0$ if and only if $\mathbf{B}\boldsymbol{\beta}=\mathbf{0}$. The coefficient vectors of the functions in $V_h\cap H^1_0(\Nb^\ell(K))$ that satisfy the constraint $\calI_H v=0$ of~\eqref{eq:mincorr} therefore form the linear subspace
\begin{equation*}
	\mathcal{A}^\ell:=\ker\mathbf{B}=\big\{\boldsymbol{\beta}\in\R^{n_0}\;:\;\mathbf{B}\boldsymbol{\beta}=\mathbf{0}\big\},
\end{equation*}
which we call the admissible set.}

\paragraph{{Discrete correction}}
{Let $\mathbf{S}\in\R^{n_0\times n_0}$ be the block of the stiffness matrix~$S_h^\ell(\bp_k)$ associated with the interior nodes, and let $\mathbf{r}\in\R^{n_0}$ be the vector with the entries $r_i:=\int_K a(x,\bp_k)\nabla\Lambda_j\cdot\nabla\lambda_{z_i}\dx$, $i=1,\dots,n_0$, such that $\mathbf{r}^\top\boldsymbol{\beta}=\int_K a(x,\bp_k)\nabla\Lambda_j\cdot\nabla v_{\boldsymbol{\beta}}\dx$ is the second term in~\eqref{eq:energyLoss}. We set
\begin{equation*}
	\mathcal{E}(\boldsymbol{\beta}):=\tfrac12\,\boldsymbol{\beta}^\top\mathbf{S}\boldsymbol{\beta}-\mathbf{r}^\top\boldsymbol{\beta},
	\qquad
	\|\boldsymbol{\beta}\|_{\mathbf{S}}:=(\boldsymbol{\beta}^\top\mathbf{S}\boldsymbol{\beta})^{1/2}=\frakA(v_{\boldsymbol{\beta}},v_{\boldsymbol{\beta}})^{1/2}
	\quad\text{for }\boldsymbol{\beta}\in\R^{n_0} .
\end{equation*}
Since $a(\cdot,\bp_k)$ is bounded from below by a positive constant and $v_{\boldsymbol{\beta}}$ vanishes on $\partial\Nb^\ell(K)$, $\|\cdot\|_{\mathbf{S}}$ is a norm on~$\R^{n_0}$. Hence, $\mathcal{E}$ is strictly convex on~$\mathcal{A}^\ell$ and has a unique minimizer
\begin{equation}\label{eq:discreteCorrector}
	\boldsymbol{\alpha}^\star:=\argmin_{\boldsymbol{\beta}\in\mathcal{A}^\ell}\;\mathcal{E}(\boldsymbol{\beta}),
\end{equation}
which is the vector of the values at the interior nodes of the solution of~\eqref{eq:mincorr} in the fine finite element space $V_h\cap H^1_0(\Nb^\ell(K))$.}

\paragraph{{Projection}}
{The values of the localized network output at the interior nodes form an arbitrary vector of~$\R^{n_0}$ and, in general, violate the interpolation condition. We therefore seek a linear map that maps every vector of~$\R^{n_0}$ into~$\mathcal{A}^\ell$, so that the condition holds by construction, and that does not exclude any admissible vector, so that $\boldsymbol{\alpha}^\star$ can still be represented. Since $\mathcal{A}^\ell=\ker\mathbf{B}$, such a map can be built from an orthonormal basis of the row space of~$\mathbf{B}$, which is the orthogonal complement of its kernel.}

\begin{lemma}[{Projection onto the admissible set}]\label{lem:projection}
{Let $r:=\rank\mathbf{B}$, let $\mathbf{Q}\in\R^{n_0\times r}$ have orthonormal columns with $\operatorname{range}\mathbf{Q}=\operatorname{range}\mathbf{B}^\top$, and define
\begin{equation*}
	P^\ell:=\Id-\mathbf{Q}\mathbf{Q}^\top\in\R^{n_0\times n_0}.
\end{equation*}
Then the following assertions hold.
\begin{enumerate}
	\item[(a)] $(P^\ell)^\top=P^\ell$ and $(P^\ell)^2=P^\ell$.
	\item[(b)] $P^\ell\boldsymbol{\beta}\in\mathcal{A}^\ell$ for all $\boldsymbol{\beta}\in\R^{n_0}$, and $P^\ell\boldsymbol{\beta}=\boldsymbol{\beta}$ for all $\boldsymbol{\beta}\in\mathcal{A}^\ell$. In particular, $P^\ell(\R^{n_0})=\mathcal{A}^\ell$.
	\item[(c)] For all $\boldsymbol{\beta}\in\R^{n_0}$, it holds $|\boldsymbol{\beta}-P^\ell\boldsymbol{\beta}|=\min_{\boldsymbol{\gamma}\in\mathcal{A}^\ell}|\boldsymbol{\beta}-\boldsymbol{\gamma}|$, where $|\cdot|$ denotes the Euclidean norm.
	\item[(d)] $P^\ell=\Id-\mathbf{B}^+\mathbf{B}$, where $\mathbf{B}^+$ denotes the Moore--Penrose inverse of~$\mathbf{B}$. In particular, $P^\ell$ does not depend on the choice of~$\mathbf{Q}$.
\end{enumerate}}
\end{lemma}
\begin{proof}
{See Appendix~\ref{app:proofProjection}.}
\end{proof}

\paragraph{{Energy of projected vectors}}
{Lemma~\ref{lem:projection} guarantees that projected vectors are admissible. It remains to show that minimizing the energy over projected vectors yields the discrete correction~$\boldsymbol{\alpha}^\star$, which is the content of the following result.}

\begin{proposition}[{Energy identity}]\label{prop:energyIdentity}
{For all $\boldsymbol{\beta}\in\R^{n_0}$, it holds
\begin{equation}\label{eq:energyIdentity}
	\mathcal{E}(P^\ell\boldsymbol{\beta})-\mathcal{E}(\boldsymbol{\alpha}^\star)=\tfrac12\,\big\|P^\ell\boldsymbol{\beta}-\boldsymbol{\alpha}^\star\big\|_{\mathbf{S}}^2 .
\end{equation}
In particular, $\min_{\boldsymbol{\beta}\in\R^{n_0}}\mathcal{E}(P^\ell\boldsymbol{\beta})=\mathcal{E}(\boldsymbol{\alpha}^\star)$, and $\boldsymbol{\beta}\in\R^{n_0}$ is a minimizer if and only if $P^\ell\boldsymbol{\beta}=\boldsymbol{\alpha}^\star$.}
\end{proposition}
\begin{proof}
{See Appendix~\ref{app:proofEnergy}.}
\end{proof}

{The excess energy $\mathcal{E}(P^\ell\boldsymbol{\beta})-\mathcal{E}(\boldsymbol{\alpha}^\star)$ is an absolute quantity that depends on the scaling of the coefficient. For reporting the accuracy of the network corrections, a relative error is more meaningful. It can be expressed by energy values only, as stated next.}

\begin{corollary}[{Loss and error}]\label{cor:lossError}
{It holds $\mathcal{E}(\boldsymbol{\alpha}^\star)=-\tfrac12\|\boldsymbol{\alpha}^\star\|_{\mathbf{S}}^2$. If $\boldsymbol{\alpha}^\star\neq\mathbf{0}$, then for all $\boldsymbol{\beta}\in\R^{n_0}$, we have
\begin{equation*}
	\frac{\|P^\ell\boldsymbol{\beta}-\boldsymbol{\alpha}^\star\|_{\mathbf{S}}}{\|\boldsymbol{\alpha}^\star\|_{\mathbf{S}}}
	=\bigg(\frac{\mathcal{E}(P^\ell\boldsymbol{\beta})-\mathcal{E}(\boldsymbol{\alpha}^\star)}{|\mathcal{E}(\boldsymbol{\alpha}^\star)|}\bigg)^{1/2}.
\end{equation*}}
\end{corollary}
\begin{proof}
{See Appendix~\ref{app:proofLossError}.}
\end{proof}

{Summed over the components and the samples, the relative error of the network corrections in the energy norm can thus be computed from the loss and the minimal energy. We use this in Section~\ref{sec:examples} to monitor the error of the network corrections during the training.}

\subsection{{Neural correction problem and DR-LOD basis}}\label{subsec:problem}
{For each element~$K\in\calT_H$, we approximate the corrections of the basis functions~$\Lambda_j$, $j\in J(K)$ (see Definition~\ref{def:lossFunction}), by one network function $v_{\btheta}^{\ell,K} := \NN$ with the $2^d$ output components $(v_{\btheta,j}^{\ell,K})_{j\in J(K)}$ and trainable network parameters $\btheta \in \R^{N_\Theta}$ according to Definition~\ref{def:RealFunction}, localized according to Definition~\ref{def:locNF}. Its input consists of the spatial coordinates and the vector $\bp_K\in\R^{N_p^K}$ of those parameters that affect the coefficient on the patch~$\Nb^\ell(K)$, i.e., $\bp_K=\bp$ if the coefficient on the patch depends on all parameters. Using the approximation~\eqref{eq:nodalInterp}, the input dataset reads $\widecheck{\bx}:= \bz^\ell \times \bp_K \in \R^{(d +N_p^K) \times n}$ with $n= N^\ell_h N_\mathrm{s}$, and we denote by $\boldsymbol{\beta}_{k,j}(\btheta)\in\R^{n_0}$ the values of the localized output of the component~$j$ for the sample~$\bp_k$ at the interior nodes of the patch, cf.~Section~\ref{subsec:exact}. Since $\mathbf{B}$ does not depend on~$\bp$, the same projection~$P^\ell$ applies to all samples. Finally, $\mathcal{E}_{k,j}$ denotes the functional~$\mathcal{E}$ of Section~\ref{subsec:exact} for the sample~$\bp_k$ and the basis function~$\Lambda_j$, i.e., the contribution of~$\bp_k$ to the discrete energy~\eqref{eq:approxLossEnergy}, restricted to the patch.}

\begin{problem}[Neural correction problem]\label{eq:optimizationProblem} {Let $\ell \geq 1$, $K \in \calT_H$, the samples $\bp_1,\dots,\bp_{N_\mathrm{s}}$, and the localized network function~$v_{\btheta}^{\ell,K}$ of the element~$K$ with the output components $v_{\btheta,j}^{\ell,K}$, $j\in J(K)$, be given, see the beginning of this subsection. Find $\btheta^\star\in\R^{N_\Theta}$ such that
\begin{equation}\label{eq:optimizationEq}
	\btheta^\star\in\argmin_{\btheta\in\R^{N_\Theta}}\;\widehat{\mathcal{L}}(\btheta),
	\qquad
	\widehat{\mathcal{L}}(\btheta):=\sum_{j\in J(K)}\sum_{k=1}^{N_\mathrm{s}}\mathcal{E}_{k,j}\big(P^\ell\boldsymbol{\beta}_{k,j}(\btheta)\big),
\end{equation}
where $\boldsymbol{\beta}_{k,j}(\btheta)\in\R^{n_0}$ contains the values of the output component~$v_{\btheta,j}^{\ell,K}(\cdot,\bp_k)$ at the interior nodes of the patch~$\Nb^\ell(K)$.}
\end{problem}

{The loss~$\widehat{\mathcal{L}}$ is the discrete energy~\eqref{eq:approxLossEnergy} of the projected network output, restricted to the patch and summed over the $2^d$ components. After the training, the network provides the corrections for an arbitrary parameter~$\bp\in\R^{N_p}$, not only for the samples: for $j\in J(K)$, let $\boldsymbol{\beta}_j(\bp)\in\R^{n_0}$ denote the values of the output component~$j$ of the trained (localized) network at the interior nodes of the patch for the input~$\bp$. The DR-LOD correction of~$\Lambda_j$ on the patch is then the finite element function with the projected values $P^\ell\boldsymbol{\beta}_j(\bp)$ at the interior nodes $z_1,\dots,z_{n_0}$,
\begin{equation*}
	v_{\btheta^\star,j}^{\ell,K}(x,\bp):=\sum_{i=1}^{n_0}\big(P^\ell\boldsymbol{\beta}_j(\bp)\big)_i\,\lambda_{z_i}(x),
	\qquad x\in\Nb^\ell(K),
\end{equation*}
where $\lambda_{z_i}$ denotes the nodal basis function of~$V_h$ associated with~$z_i$. By construction, the DR-LOD correction $v_{\btheta^\star,j}^{\ell,K}(\cdot,\bp)$ belongs to $V_h\cap H^1_0(\Nb^\ell(K))$ and satisfies the interpolation condition $\calI_H v_{\btheta^\star,j}^{\ell,K}(\cdot,\bp)=0$, cf.~Section~\ref{subsec:exact}.}

{Summing the energy identity~\eqref{eq:energyIdentity} over the components and the samples, with the discrete corrections~$\boldsymbol{\alpha}^\star_{k,j}$ of~\eqref{eq:discreteCorrector} and the corresponding matrices~$\mathbf{S}_k$, yields
\begin{equation*}
	\widehat{\mathcal{L}}(\btheta)
	=\widehat{\mathcal{L}}^\star
	+\frac12\sum_{j\in J(K)}\sum_{k=1}^{N_\mathrm{s}}\big\|P^\ell\boldsymbol{\beta}_{k,j}(\btheta)-\boldsymbol{\alpha}^\star_{k,j}\big\|^2_{\mathbf{S}_k},
	\qquad
	\widehat{\mathcal{L}}^\star:=\sum_{j\in J(K)}\sum_{k=1}^{N_\mathrm{s}}\mathcal{E}_{k,j}(\boldsymbol{\alpha}^\star_{k,j}) .
\end{equation*}
Since $\widehat{\mathcal{L}}^\star$ does not depend on~$\btheta$, minimizing the loss is equivalent to minimizing the sum of the squared errors of the network corrections in the energy norm over all components and samples. Hence, minimizing the loss amounts to computing the best approximation of the discrete corrections $\balpha_{k,j}^\star$ by the network in this norm without requiring $\balpha_{k,j}^\star$.}

In a post-processing step, i.e., after training the ANN by solving~\eqref{eq:optimizationEq}, we compute the corrected multiscale basis $\{\widehat{\Lambda}_j^\ell\}_{j=1}^{N_H}$ that spans the modified multiscale space $ \widehat{V}_H^\ell$. That is, each $\widehat{\Lambda}_j^\ell$ is associated to the corresponding basis function $\Lambda_j $ of $ V_H$ and related to the coefficient function~$a\colon\Omega \times \R^{N_p} \to \R$ induced by $\bp$. More precisely, it is obtained similarly to~\eqref{eq:corrBF} by 
\begin{equation}\label{eq:annBasis}
	\big(\widehat{\Lambda}_j^\ell(\bp)\big)(x) = \Lambda_j(x) - \sum_{{K}\in\calT_H\colon\operatorname{supp}(\Lambda_j)\cap {K}\neq\emptyset}v_{\btheta,j}^{\ell,{K}}(x,\bp){,}
\end{equation}
{where $v_{\btheta^\star,j}^{\ell,K}$ denotes the component~$j$ of the trained network of the element~$K$, see Problem~\ref{eq:optimizationProblem}.}
By computing/predicting the multiscale basis functions corresponding to all $\Lambda_j \in V_H$, $j=1,\dots,N_H$, solving equation~\eqref{eq:fd_heat} follows the lines of the classical LOD approach, cf.~\cite{MalP20}.
That is, for an elliptic model problem and given some parameter $\bp$, we compute the Galerkin solution in the space $\widehat V_H^\ell$. This requires the computation of appropriate system matrices making use of the corrected basis functions. For time-dependent or uncertain problems, such matrices need to be updated based on the choice of $\bp$ during simulations. To compute such matrices more efficiently, a Petrov--Galerkin variant of the LOD may be used, see also~\cite{EngHMP19}.

\subsection{{Training procedure and implementation}}\label{subsec:trainingDetails}
{We summarize the complete procedure, see also Algorithm~\ref{alg:lodann}. For each coarse element~$K$, one network is trained according to Problem~\ref{eq:optimizationProblem}. The networks of different elements are independent of each other and can be trained in parallel.}

{Finally, we normalize the inputs and the output of the network by fixed (non-trainable) affine layers. The spatial input is mapped from the respective patch onto $[-1,1]^d$ and the parameters from their range onto $[-1,1]$. The output component~$j\in J(K)$ is multiplied by a fixed scale~$s_{K,j}$, since the corrections are typically much smaller than the output of an initialized network. The scale is determined from the energy~\eqref{eq:optimizationEq} itself: denoting by $\boldsymbol{\beta}^0_{k,j}$ the values of the localized output of the initialized network (without scaling) at the interior nodes, $s_{K,j}$ is the exact line search of the energy along the projected output,
\begin{equation*}
	s_{K,j} := \argmin_{s\in\R}\sum_{k\in\mathcal{M}} \mathcal{E}_{k,j}\big(s\,P^\ell\boldsymbol{\beta}^0_{k,j}\big)
	= \frac{\sum_{k\in\mathcal{M}} \mathbf{r}_{k,j}^\top P^\ell\boldsymbol{\beta}^0_{k,j}}{\sum_{k\in\mathcal{M}} \big\|P^\ell\boldsymbol{\beta}^0_{k,j}\big\|_{\mathbf{S}_k}^2},
\end{equation*}
where $\mathbf{r}_{k,j}$ denotes the right-hand side of~$\mathcal{E}_{k,j}$ and $\mathcal{M}\subset\{1,\dots,N_\mathrm{s}\}$ a small subset of the samples ($16$ equidistant indices in our experiments). This requires only matrix-vector products with the matrices of the training. Since the scaling only multiplies the last layer of the network, it does not change the set of network functions, and hence neither the minimizers of~\eqref{eq:optimizationEq}, but only the conditioning of the optimization. The matrix $\mathbf{Q}$ of Lemma~\ref{lem:projection} is computed once per patch from a singular value decomposition of~$\mathbf{B}$ as the right singular vectors to the non-zero singular values, where the rank is determined by a relative threshold on the singular values. The projection~$P^\ell$ is implemented as a fixed (non-trainable) layer of the network; applying it to the values at the interior nodes costs $\mathcal{O}(n_0 r)$ operations. In our experiments, the normalization substantially accelerates the training.}

{The trainable parameters are determined by gradient-based optimization, where the gradients with respect to~$\btheta$ are computed by backpropagation, a special case of the reverse mode of automatic differentiation~\cite{rumelhart1986learning,griewank2003mathematical,baydin2018automatic}. We combine the adaptive moment estimation (Adam)~\cite{kingma2015adam}, a stochastic gradient method~\cite{bottou2018optimization}, with the limited-memory Broyden--Fletcher--Goldfarb--Shanno method (L-BFGS)~\cite{liu1989limited}, which is a common combination for the training of PINNs. More precisely, we first train with Adam and an exponentially decaying learning rate. Adam is stopped as soon as the loss has not decreased by more than a relative tolerance of~$10^{-4}$ within $2000$ consecutive epochs (or after a maximal number of epochs). Afterwards, the network is fine-tuned with L-BFGS (memory size~$50$) in double precision for a fixed number of iterations. The fine-tuning considerably improves the accuracy of the trained corrections compared to a training with Adam only.}

{If the number of samples~$N_\mathrm{s}$ is large, each Adam epoch uses a randomly chosen mini-batch of coefficients, and L-BFGS uses a fixed subset of the samples. Apart from the coefficients for the parameter samples, no data is needed for the training.} 

\begin{algorithm}[t]
{\textbf{Offline phase} (training)
\begin{enumerate}[label=\arabic*.,leftmargin=*,itemsep=2pt]
	\item Choose parameter samples $\bp_1,\dots,\bp_{N_\mathrm{s}}$.
	\item For each coarse element $K\in\calT_H$ (in parallel):
	\begin{enumerate}[label=(\alph*),itemsep=1pt]
		\item assemble the stiffness matrices $S_h^\ell(\bp_k)$, $k=1,\dots,N_\mathrm{s}$, on the patch $\Nb^\ell(K)$ and the right-hand sides of~\eqref{eq:optimizationEq};
		\item compute the projection~$P^\ell$ from the local interpolation matrix (Lemma~\ref{lem:projection}) and the output scale (one local solve);
		\item set up the network with the input normalization, the output scaling, the localization~$g$, and the projection~$P^\ell$ as fixed layers;
		\item minimize the energy~\eqref{eq:optimizationEq} with Adam until it stagnates (or max iterations is reached);
		\item fine-tune with L-BFGS in double precision.
	\end{enumerate}
\end{enumerate}
\textbf{Online phase} (for each time step $m$ with parameter~$\bp^m$)
\begin{enumerate}[label=\arabic*.,leftmargin=*,itemsep=2pt,start=3]
	\item Evaluate the networks at the fine nodes of the patches for~$\bp^m$ and apply~$P^\ell$.
	\item Assemble the basis functions $\widehat\Lambda_j^\ell(\bp^m)$ by~\eqref{eq:annBasis} and the coarse system matrices.
	\item Solve the coarse backward Euler step for $\widehat u_H^m\in\widehat V_H^\ell$.
\end{enumerate}}
\caption{{Offline and online phase of the DR-LOD method.}}\label{alg:lodann}
\end{algorithm}

%----------------------------------------------------------
% Examples
%----------------------------------------------------------
\section{Numerical Examples}\label{sec:examples}

In order to apply our method to the parabolic equation~\eqref{eq:heat_pde}, we consider the variational formulation~\eqref{eq:heat_var} and discretize in space and time. For the spatial discretization, we first employ for a given~$\bp$ the classical LOD space~$\widetilde{V}_H^\ell$ as given in Definition~\ref{def:locCorrOp} combined with a backward Euler method {(see also\eqref{eq:fd_heat_lod})}. The corresponding approximate solution in the $m$th time step is denoted by~$\widetilde{u}_H^m$ and {and serves as a benchmark for the DR-LOD method.} 
Second, we compute a {DR-LOD} space~$\widehat{V}_H^\ell$ for the spatial approximation based on the derivations in Section~\ref{sec:ann}, again based on a backward Euler method in time. The corresponding approximate solution in the $m$th time step is denoted by~$\widehat{u}_H^m$. 
Note that the time step for both approaches is given by~$\tau={T}/{M}$, where~$M$ is the number of time steps. {Since the classical LOD method approximates the solution of the parabolic model problem accurately, if the solution of the {DR-LOD} approach is close to the LOD solution, we have a good approximation of the exact solution as well.}
All numerical examples have been conducted using the code available at \url{https://github.com/melasmi8/drlod.git}. \\
{In both examples, we compare the DR-LOD solution $\widehat{u}_H^m$ and the LOD solution $\widetilde{u}_H^m$ with the fine finite element solution~$u_h^m$ of~\eqref{eq:fd_heat_fem}, which serves as the reference solution. Further, we consider the classical finite element method on the coarse mesh (coarse FEM) and an LOD method whose basis is computed once for the initial coefficient and then kept fixed for all time steps (frozen LOD). The latter shows the error that is caused if the multiscale basis is not updated when the coefficient changes. All errors are relative errors in the Euclidean norm of the nodal values on the fine mesh, i.e., a discrete $L^2$-norm, which is computed at every time step. Similarly, the error of the DR-LOD basis functions is the relative error of their nodal values compared to the LOD basis functions, averaged over all basis functions and time steps.

The networks are fully connected with width~$128$ and depth~$8$. They are trained as described in Section~\ref{subsec:trainingDetails}. Adam is used with an exponentially decaying learning rate (initial learning rate $10^{-3}$, decay rate $0.9$ every $1000$ epochs) for at most $30\,000$ epochs.

The relevant cost of the approach is the cost of the offline phase, i.e., the training. Once the networks are trained, the DR-LOD basis for a new coefficient is obtained by a forward pass of the networks and the projection~$P^\ell$, without solving any local problem.}

\subsection{{Example 1: damage front}}\label{subsec:ex1}
{In the first example, we consider a heterogeneous material with a growing damage front. The coefficient is given by
\begin{equation*}
	a(x,p) = a_\mathrm{base}(x)\,\big(1-(1-\kappa)\,D(x,p)\big)
\end{equation*}
\begin{equation*}
 \text{with } 	D(x,p) = \Big(1+\exp\Big(-\frac{\eta\,p-\phi(x)}{\delta}\Big)\Big)^{-1},
	\qquad p\in[0,1],
\end{equation*}
and $a_\mathrm{base}$ is piecewise constant on a uniform grid with $\varepsilon=1/36$ with values that are log-uniformly distributed in $[0.3,1]$, $\kappa=0.1$, $\eta=0.7$, and $\delta=0.015$. The function $\phi$ is a smooth random field with correlation length of about $0.1$ (a superposition of $60$ random cosine modes), which is transformed such that its values are uniformly distributed in $[0,1]$. The function $D\in(0,1)$ indicates the damaged regions, in which the conductivity is reduced by the factor~$\kappa$. The damaged area fraction is about $\eta\,p$, i.e., as $p$ increases, damaged regions nucleate, grow, and merge, see Figure~\ref{fig:ex1coeff}. Hence, the coefficient varies on the scale~$\varepsilon$ with values in about $[0.03,1]$, and its dependence on~$p$ is strongly nonlinear.

We choose $f\equiv1$, $u_0\equiv0$, $T=1$, $H=1/12$, $h=1/36$, $\ell=2$, and the time step size $\tau=1/24$. The parameter, and thus the damage, increases in time: for the test set, we draw $25$ parameters uniformly at random in $[0,1]$, sort them increasingly, and use one of them in each time step $t^m=m\tau$, $m=0,\dots,24$. For the training, we use $N_\mathrm{s}=40$ equidistant parameters in $[0,1]$, which do not coincide with the test parameters. Adam stopped after $2900$ to $9000$ epochs (median $6500$) for the $144$ networks, and L-BFGS was run for $2000$ iterations. A typical training history is depicted in Figure~\ref{fig:ex1loss}. Besides the energy loss, it shows the relative error of the network corrections of all training samples in the energy norm, which we compute from the loss and the minimal energy by Corollary~\ref{cor:lossError}. For the element shown, the error is $1.1\,\%$ after the Adam phase, and L-BFGS reduces it to $0.5\,\%$. The training of all $144$ networks took $5.6$ hours on a workstation with $32$ CPU cores, with $12$ networks trained in parallel processes.}

\begin{figure}[htbp]
	\centering
	\begin{subfigure}[b]{0.42\textwidth}\centering
		\input{./coefficients_0.tikz}
		\caption{$p=0$ (training)}\label{fig:ex1coefficients-a}
	\end{subfigure}\hspace{0.05\textwidth}
	\begin{subfigure}[b]{0.42\textwidth}\centering
		\input{./coefficients_1.tikz}
		\caption{$p=20/39$ (training)}\label{fig:ex1coefficients-b}
	\end{subfigure}
	\par\smallskip
	\begin{subfigure}[b]{0.42\textwidth}\centering
		\input{./coefficients_2.tikz}
		\caption{$p=1$ (training)}\label{fig:ex1coefficients-c}
	\end{subfigure}\hspace{0.05\textwidth}
	\begin{subfigure}[b]{0.42\textwidth}\centering
		\input{./coefficients_3.tikz}
		\caption{$p\approx0.45$ (test)}\label{fig:ex1coefficients-d}
	\end{subfigure}
	\caption{{Coefficients of Example~1 for three training parameters and for a test parameter.}}\label{fig:ex1coeff}
\end{figure}

\begin{figure}[htbp]
	\centering
	\begin{subfigure}[b]{0.42\textwidth}\centering
		\input{./loss_0.tikz}
		\caption{energy loss}\label{fig:ex1loss-a}
	\end{subfigure}\hfill
	\begin{subfigure}[b]{0.42\textwidth}\centering
		\input{./loss_1.tikz}
		\caption{relative error of the corrections}\label{fig:ex1loss-b}
	\end{subfigure}
	\caption{{Training history of the network of the element $K=(1/2,7/12)^2$ in Example~1: energy loss~$\widehat{\mathcal{L}}$ of~\eqref{eq:optimizationEq} (a) and relative error of the network corrections in the energy norm, $((\widehat{\mathcal{L}}-\widehat{\mathcal{L}}^\star)/|\widehat{\mathcal{L}}^\star|)^{1/2}$ with the minimal value~$\widehat{\mathcal{L}}^\star$ of the loss}, see Corollary~\ref{cor:lossError} (b), during the Adam (solid) and the L-BFGS (dashed) phase.}\label{fig:ex1loss}
\end{figure}

{Figure~\ref{fig:ex1basis} compares three LOD basis functions with their DR-LOD counterparts for a test parameter. The networks capture the fine-scale oscillations of the corrected basis functions, which differ considerably from the classical finite element basis functions. The mean relative error of the DR-LOD basis functions is $0.15\,\%$ for the test parameters and $0.16\,\%$ for the training parameters, i.e., the networks generalize to parameters that have not been seen during training.}
\begin{figure}[htbp]
	\centering
	\begin{subfigure}[b]{0.32\textwidth}\centering
		\input{./basis_functions_0.tikz}
		\caption{LOD, node $(1/2,1/2)$}\label{fig:ex1basisfunctions-a}
	\end{subfigure}\hfill
	\begin{subfigure}[b]{0.32\textwidth}\centering
		\input{./basis_functions_1.tikz}
		\caption{LOD, node $(1/12,1/12)$}\label{fig:ex1basisfunctions-b}
	\end{subfigure}\hfill
	\begin{subfigure}[b]{0.32\textwidth}\centering
		\input{./basis_functions_2.tikz}
		\caption{LOD, node $(1/2,1/12)$}\label{fig:ex1basisfunctions-c}
	\end{subfigure}
	\par\smallskip
	\begin{subfigure}[b]{0.32\textwidth}\centering
		\input{./basis_functions_3.tikz}
		\caption{DR-LOD, node $(1/2,1/2)$}\label{fig:ex1basisfunctions-d}
	\end{subfigure}\hfill
	\begin{subfigure}[b]{0.32\textwidth}\centering
		\input{./basis_functions_4.tikz}
		\caption{DR-LOD, node $(1/12,1/12)$}\label{fig:ex1basisfunctions-e}
	\end{subfigure}\hfill
	\begin{subfigure}[b]{0.32\textwidth}\centering
		\input{./basis_functions_5.tikz}
		\caption{DR-LOD, node $(1/2,1/12)$}\label{fig:ex1basisfunctions-f}
	\end{subfigure}
	\par\smallskip
	\begin{subfigure}[b]{0.32\textwidth}\centering
		\input{./basis_functions_6.tikz}
		\caption{difference, node $(1/2,1/2)$}\label{fig:ex1basisfunctions-g}
	\end{subfigure}\hfill
	\begin{subfigure}[b]{0.32\textwidth}\centering
		\input{./basis_functions_7.tikz}
		\caption{difference, node $(1/12,1/12)$}\label{fig:ex1basisfunctions-h}
	\end{subfigure}\hfill
	\begin{subfigure}[b]{0.32\textwidth}\centering
		\input{./basis_functions_8.tikz}
		\caption{difference, node $(1/2,1/12)$}\label{fig:ex1basisfunctions-i}
	\end{subfigure}
	\caption{{LOD basis functions, DR-LOD basis functions, and their difference in Example~1 for three nodes and the test parameter $p\approx0.45$.}}\label{fig:ex1basis}
\end{figure}
{In Figure~\ref{fig:ex1sol}, we depict the coefficient, the LOD solution, and the DR-LOD solution for the test set at the times $t=1/24$, $t=1/2$, and $t=1$. The cross sections in Figure~\ref{fig:ex1solCS} and the pointwise errors in Figure~\ref{fig:ex1err} show that the DR-LOD solution is visually indistinguishable from the LOD solution. In contrast, the coarse FEM and the frozen LOD solution deviate considerably from the reference solution. This is quantified in Table~\ref{tab:ex1}. The DR-LOD solution deviates from the LOD solution by only $0.06\,\%$ at the final time, such that its error compared to the reference solution ($0.97\,\%$) is the error of the LOD itself ($0.97\,\%$) for the chosen $H$ and~$\ell$. The frozen LOD, i.e., the LOD without updates of the basis, has an error of $8.8\,\%$, about nine times the error of the LOD, and the coarse FEM an error of $12.2\,\%$, which illustrates that the multiscale basis has to be updated whenever the coefficient changes.}
\begin{table}[htbp]
	\centering
	{\begin{tabular}{l c c c c}
		\hline
		& \multicolumn{2}{c}{training parameters} & \multicolumn{2}{c}{test parameters} \\
		method & final time & mean over time & final time & mean over time \\
		\hline\hline
		coarse FEM & $12.19\,\%$ & $9.36\,\%$ & $12.20\,\%$ & $9.00\,\%$ \\
		frozen LOD & $8.47\,\%$ & $6.07\,\%$ & $8.81\,\%$ & $5.87\,\%$ \\
		LOD & $0.94\,\%$ & $1.01\,\%$ & $0.97\,\%$ & $0.98\,\%$ \\
		DR-LOD & $0.95\,\%$ & $1.01\,\%$ & $0.97\,\%$ & $0.98\,\%$ \\
		\hline
		DR-LOD vs.\ LOD & $0.07\,\%$ & $0.04\,\%$ & $0.06\,\%$ & $0.04\,\%$ \\
		\hline
	\end{tabular}}
	\caption{{Relative errors compared to the fine finite element solution in Example~1 at the final time and averaged over all time steps, and relative difference between the DR-LOD and the LOD solution (last row). For the training parameters, the $40$ training coefficients are used in consecutive time steps with $\tau=1/39$.}}\label{tab:ex1}
\end{table}

\begin{figure}[htbp]
	\centering
	\begin{subfigure}[b]{0.32\textwidth}\centering
		\input{./solutions_0.tikz}
		\caption{coefficient, $t=1/24$}\label{fig:ex1solutions-a}
	\end{subfigure}\hfill
	\begin{subfigure}[b]{0.32\textwidth}\centering
		\input{./solutions_1.tikz}
		\caption{coefficient, $t=1/2$}\label{fig:ex1solutions-b}
	\end{subfigure}\hfill
	\begin{subfigure}[b]{0.32\textwidth}\centering
		\input{./solutions_2.tikz}
		\caption{coefficient, $t=1$}\label{fig:ex1solutions-c}
	\end{subfigure}
	\par\smallskip
	\begin{subfigure}[b]{0.32\textwidth}\centering
		\input{./solutions_6.tikz}
		\caption{LOD, $t=1/24$}\label{fig:ex1solutions-g}
	\end{subfigure}\hfill
	\begin{subfigure}[b]{0.32\textwidth}\centering
		\input{./solutions_7.tikz}
		\caption{LOD, $t=1/2$}\label{fig:ex1solutions-h}
	\end{subfigure}\hfill
	\begin{subfigure}[b]{0.32\textwidth}\centering
		\input{./solutions_8.tikz}
		\caption{LOD, $t=1$}\label{fig:ex1solutions-i}
	\end{subfigure}
	\par\smallskip
	\begin{subfigure}[b]{0.32\textwidth}\centering
		\input{./solutions_9.tikz}
		\caption{DR-LOD, $t=1/24$}\label{fig:ex1solutions-j}
	\end{subfigure}\hfill
	\begin{subfigure}[b]{0.32\textwidth}\centering
		\input{./solutions_10.tikz}
		\caption{DR-LOD, $t=1/2$}\label{fig:ex1solutions-k}
	\end{subfigure}\hfill
	\begin{subfigure}[b]{0.32\textwidth}\centering
		\input{./solutions_11.tikz}
		\caption{DR-LOD, $t=1$}\label{fig:ex1solutions-l}
	\end{subfigure}
	\caption{{Coefficient, LOD solution, and DR-LOD solution for the test set of Example~1.}}\label{fig:ex1sol}
\end{figure}

\begin{figure}[htbp]
	\centering
	\begin{subfigure}[b]{0.31\textwidth}\centering
		\input{./solution_cross_sections_0.tikz}
		\caption{$t=1/24$}\label{fig:ex1solutioncrosssections-a}
	\end{subfigure}\hfill
	\begin{subfigure}[b]{0.31\textwidth}\centering
		\input{./solution_cross_sections_1.tikz}
		\caption{$t=1/2$}\label{fig:ex1solutioncrosssections-b}
	\end{subfigure}\hfill
	\begin{subfigure}[b]{0.31\textwidth}\centering
		\input{./solution_cross_sections_2.tikz}
		\caption{$t=1$}\label{fig:ex1solutioncrosssections-c}
	\end{subfigure}
	\par\smallskip
	\input{./solution_cross_sections_legend.tikz}
	\caption{{Cross sections at $y=1/2$ of the solutions of all methods for the test set of Example~1 at the times $t=1/24$, $t=1/2$, and $t=1$.}}\label{fig:ex1solCS}
\end{figure}

\begin{figure}[htbp]
	\centering
	\begin{subfigure}[b]{0.32\textwidth}\centering
		\input{./pointwise_errors_0.tikz}
		\caption{LOD error, $t=1/24$}\label{fig:ex1pointwiseerrors-a}
	\end{subfigure}\hfill
	\begin{subfigure}[b]{0.32\textwidth}\centering
		\input{./pointwise_errors_1.tikz}
		\caption{LOD error, $t=1/2$}\label{fig:ex1pointwiseerrors-b}
	\end{subfigure}\hfill
	\begin{subfigure}[b]{0.32\textwidth}\centering
		\input{./pointwise_errors_2.tikz}
		\caption{LOD error, $t=1$}\label{fig:ex1pointwiseerrors-c}
	\end{subfigure}
	\par\smallskip
	\begin{subfigure}[b]{0.32\textwidth}\centering
		\input{./pointwise_errors_3.tikz}
		\caption{DR-LOD error, $t=1/24$}\label{fig:ex1pointwiseerrors-d}
	\end{subfigure}\hfill
	\begin{subfigure}[b]{0.32\textwidth}\centering
		\input{./pointwise_errors_4.tikz}
		\caption{DR-LOD error, $t=1/2$}\label{fig:ex1pointwiseerrors-e}
	\end{subfigure}\hfill
	\begin{subfigure}[b]{0.32\textwidth}\centering
		\input{./pointwise_errors_5.tikz}
		\caption{DR-LOD error, $t=1$}\label{fig:ex1pointwiseerrors-f}
	\end{subfigure}
	\par\smallskip
	\begin{subfigure}[b]{0.32\textwidth}\centering
		\input{./pointwise_errors_6.tikz}
		\caption{DR-LOD $-$ LOD, $t=1/24$}\label{fig:ex1pointwiseerrors-g}
	\end{subfigure}\hfill
	\begin{subfigure}[b]{0.32\textwidth}\centering
		\input{./pointwise_errors_7.tikz}
		\caption{DR-LOD $-$ LOD, $t=1/2$}\label{fig:ex1pointwiseerrors-h}
	\end{subfigure}\hfill
	\begin{subfigure}[b]{0.32\textwidth}\centering
		\input{./pointwise_errors_8.tikz}
		\caption{DR-LOD $-$ LOD, $t=1$}\label{fig:ex1pointwiseerrors-i}
	\end{subfigure}
	\caption{{Absolute pointwise errors of the LOD solution and of the DR-LOD solution compared to the fine finite element solution, and absolute pointwise difference between the DR-LOD and the LOD solution for the test set of Example~1.}}\label{fig:ex1err}
\end{figure}

{The size of the networks is a trade-off between accuracy and training cost. A study of networks of different widths and depths in~\ref{app:size}, carried out for this example on the coarser mesh $H=1/6$ with $\ell=1$, shows that the error of the network corrections decreases at an observed rate of about $N_\Theta^{-1/2}$ in the number of weights~$N_\Theta$ when the width is increased, i.e., about first order in the width, and even faster when the depth is increased. Networks of the size used here (width~$128$, depth~$8$) reach a mean relative error of the corrections of $1.0\,\%$ in the energy norm, whereas considerably smaller networks are clearly less accurate, e.g., $1.9\,\%$ for width~$64$ and depth~$4$ and $35\,\%$ for width~$16$ and depth~$2$, at a training time per network that is lower by up to a factor of about~$9$.

For a single scalar parameter as in this example, an alternative to the networks would be to precompute the LOD basis for a set of parameters and to interpolate it in~$p$. For instance, linear interpolation between the LOD bases for $40$ equidistant parameters leads to a mean relative error of the basis functions of $0.15\,\%$ for the test parameters, which is the same as the accuracy of the networks. However, the number of required local problems of such an interpolation grows exponentially with the number of parameters. This is the setting of the second example.}

\subsection{{Example 2: battery cell with multiple parameters}}\label{subsec:ex2}

For the second example, we consider the simplified model of~\cite{Ves22} in the context of battery simulations. It reads 
\begin{equation}\label{eq:heat_eq_battery}
	\rho\, \partial_t c_\mathrm{p} u-\Div (\overline\lambda \nabla u) = f \quad \text{ in } \Omega\times (0,T] {, \quad T>0}
\end{equation}
with homogeneous boundary conditions on $\partial \Omega$, where $\overline \lambda$ denotes the heat conductivity and the density~$\rho$ as well as the specific heat capacity $c_\mathrm{p}$ are assumed to be constant in time. The model is equivalent to~\eqref{eq:heat_pde} with $a={\overline\lambda}/{\rho c_\mathrm{p}}$ with a scaled right-hand side. {The domain $\Omega=(0,1)^2$ consists of horizontal layers, namely the cathode collector (CC), the active material (AM), the anode collector (AC), and again the AM, which are repeated four times and closed by a CC layer. The collector layers are four times thinner than the AM layers, such that there are eight AM layers, see Figure~\ref{fig:ex2coeff}. The material properties are given in Table~\ref{tab:material_properties}. In the computations, the volumetric heat capacity $\rho c_\mathrm{p}$ is incorporated as a weight in the mass matrices of the backward Euler method~\eqref{eq:fd_heat}, and the LOD corrections are computed for the conductivity~$\overline\lambda$.} 
\begin{table}[h!]
  \centering
	\begin{tabular}{c c c c} 
		\hline
		& $\rho$ & $c_\mathrm{p}$ & $\overline\lambda$ \\ [0.5ex] 
		\hline\hline
		AC & 8710.2 & 384.65 & {398.71} \\ 
		CC & 2706.77 & 897.8 & 236.3 \\
		AM & 2094.302 & 1010.119 & {$[2,250]$}\\ 
		\hline
	\end{tabular}
	\caption{Material properties in Example 2}\label{tab:material_properties}
\end{table}
{In contrast to the collectors, the AM is strongly heterogeneous, and its conductivity depends on the local state of the material, e.g., the local temperature or state of charge, which changes in time. To model this, we divide each of the eight AM layers into six segments in $x$-direction, which results in $48$ regions~$R_r$, $r=1,\dots,48$, each with its own state parameter $p_r\in[0,1]$. The conductivity of the AM is given by
\begin{equation*}
	\overline\lambda(x,\bp) = \lambda_\mathrm{base}(x)\, 5^{2p_r-1}\qquad\text{for } x\in R_r,\quad \bp=(p_1,\dots,p_{48})\in[0,1]^{48},
\end{equation*}
where $\lambda_\mathrm{base}=(\lambda_0\lambda_1)^{1/2}$ with two random fields $\lambda_0,\lambda_1$ that are piecewise constant on a uniform grid with $\varepsilon=1/60$ and log-uniformly distributed in $[10,50]$. That is, the conductivity of each region is scaled by a factor between $1/5$ and~$5$, and $\overline\lambda$ varies in $[2,250]$ in the AM. The coefficient thus depends on $N_p=48$ parameters, which is in line with the uncertainty assumptions mentioned at the end of Section~\ref{subsec:lod}.}
Specifically, in the offline phase prior to the simulations we might not have exact a priori knowledge on the time evolution of the coefficient, such that we have to consider a whole range of possible approximation spaces as well. This setting is where the LOD method has its bottleneck, as it cannot pre-compute the basis functions if the coefficient is not known. Moreover, even for known coefficients with rapid variations in time the computation of all necessary basis functions is computationally too demanding. In this setting the {DR-LOD} method performs optimally as we can pre-compute (in the sense of training the network) a wide range of coefficients such that obtaining the correctors for each coefficient within the parameter space can be performed fast during the online phase. 
{Such a high-dimensional parameter space is particularly challenging for the LOD method. Precomputing the LOD basis for a set of parameters and interpolating it, as discussed for the first example, is not feasible: each local correction problem depends on the parameters of all regions that intersect its patch, i.e., on $6$ to $15$ parameters. A tensor-product interpolation in these parameters would require about $10^8$ (three values per parameter) or $2\cdot10^{11}$ (five values per parameter) local correction problems. In contrast, the networks are trained with a moderate number of coefficient samples, and each network only receives the spatial coordinates and the $6$ to $15$ parameters on which its correction problem depends, cf.~Problem~\ref{eq:optimizationProblem}.

We use $H=1/6$, $\ell=1$, $h=1/60$, the source term $f=3$ in the AM and $f=0$ in the collectors, the initial condition $u_0\equiv0$, and the final time $T=10^4$ with the time step size $\tau=T/24$. Since the volumetric heat capacity $\rho c_\mathrm{p}$ is of the order $10^6$, heat spreads slowly, and a long time horizon is required for the diffusion to become relevant. For the training, we use $N_\mathrm{s}=256$ parameter vectors in $[0,1]^{48}$ drawn by Latin hypercube sampling. In each Adam epoch, a mini-batch of $16$ coefficients is used, and the subsequent $1000$ L-BFGS iterations use a fixed subset of $32$ coefficients. The total training time was about $27$ hours. For the test set, all region states evolve smoothly in time as $p_r(t)=0.5+0.45\sin(2\pi(f_rt/T+\varphi_r))$ with frequencies $f_r\in[0.5,1.5]$ and phases $\varphi_r\in[0,1]$ drawn uniformly at random, such that none of the $25$ test coefficients has been seen during training, see Figure~\ref{fig:ex2coeff} (right).}
\begin{figure}[htbp]
	\centering
	\begin{subfigure}[b]{0.42\textwidth}\centering
		\input{./ex2_coefficients_0.tikz}
		\caption{training sample $\bp_1$}\label{fig:ex2coefficients-a}
	\end{subfigure}\hspace{0.05\textwidth}
	\begin{subfigure}[b]{0.42\textwidth}\centering
		\input{./ex2_coefficients_1.tikz}
		\caption{training sample $\bp_{13}$}\label{fig:ex2coefficients-b}
	\end{subfigure}
	\par\smallskip
	\begin{subfigure}[b]{0.42\textwidth}\centering
		\input{./ex2_coefficients_2.tikz}
		\caption{training sample $\bp_{25}$}\label{fig:ex2coefficients-c}
	\end{subfigure}\hspace{0.05\textwidth}
	\begin{subfigure}[b]{0.42\textwidth}\centering
		\input{./ex2_coefficients_3.tikz}
		\caption{test set, $t=5\cdot10^3$}\label{fig:ex2coefficients-d}
	\end{subfigure}
	\caption{{Heat conductivity $\overline\lambda$ in Example~2 for three training parameter vectors (a)--(c) and for the test set at time $t=5\cdot10^3$ (d). The thin horizontal layers with constant conductivity are the collectors, while the $48$ regions of the AM have individual states.}}\label{fig:ex2coeff}
\end{figure}
{Figure~\ref{fig:ex2basis} compares three LOD basis functions with their DR-LOD counterparts for a test coefficient. The mean relative error of the DR-LOD basis functions over the test set is $3.1\,\%$, which is larger than in the first example due to the much larger parameter space. However, the resulting solutions are still very accurate. In Figure~\ref{fig:ex2sol}, we depict the coefficient and the solutions for the test set at the times $t=417$, $t=5\cdot10^3$, and $t=10^4$, and Table~\ref{tab:ex2} shows the errors at the final time and averaged over time. The DR-LOD solution deviates from the LOD solution by $0.71\,\%$ at the final time, and its error compared to the reference solution ($5.29\,\%$) is close to that of the LOD ($5.17\,\%$). The coarse FEM ($35.6\,\%$) and the frozen LOD ($22.3\,\%$) are much less accurate.

As a further comparison, we consider a classical Deep Ritz method without the LOD structure: a single network on the whole domain (width $128$, depth $8$) that minimizes in each time step the discrete energy of the fine-scale backward Euler step~\eqref{eq:fd_heat_fem}. The network is trained anew in each time step (warm started from the previous step with $500$ Adam and $500$ L-BFGS iterations, and $3000$ Adam iterations in the first step). It only learns the solution for one particular coefficient sequence and thus has no offline phase that could be reused for other coefficients. Its error at the final time is $7.1\,\%$ (mean over time $10.8\,\%$, compared to $7.6\,\%$ of DR-LOD), and the retraining requires about $167$ seconds per time step, compared to $0.14$ seconds per time step for DR-LOD, including the evaluation of the networks.}

\begin{table}[htbp]
	\centering
	{\begin{tabular}{l c c c c}
		\hline
		& \multicolumn{2}{c}{training parameters} & \multicolumn{2}{c}{test set} \\
		method & final time & mean over time & final time & mean over time \\
		\hline\hline
		coarse FEM & $34.31\,\%$ & $30.73\,\%$ & $35.61\,\%$ & $32.49\,\%$ \\
		frozen LOD & $15.71\,\%$ & $15.40\,\%$ & $22.27\,\%$ & $21.14\,\%$ \\
		Deep Ritz & -- & -- & $7.12\,\%$ & $10.82\,\%$ \\
		LOD & $4.84\,\%$ & $6.10\,\%$ & $5.17\,\%$ & $7.57\,\%$ \\
		DR-LOD & $4.87\,\%$ & $6.15\,\%$ & $5.29\,\%$ & $7.64\,\%$ \\
		\hline
		DR-LOD vs.\ LOD & $0.77\,\%$ & $0.51\,\%$ & $0.71\,\%$ & $0.70\,\%$ \\
		\hline
	\end{tabular}}
	\caption{{Relative errors compared to the fine finite element solution in Example~2 at the final time and averaged over all time steps, and relative difference between the DR-LOD and the LOD solution (last row). For the training parameters, the first $25$ training parameter vectors are used in consecutive time steps.}}\label{tab:ex2}
\end{table}

\begin{figure}[htbp]
	\centering
	\begin{subfigure}[b]{0.32\textwidth}\centering
		\input{./ex2_basis_functions_0.tikz}
		\caption{{LOD, node $(1/2,1/2)$}}\label{fig:ex2basisfunctions-a}
	\end{subfigure}\hfill
	\begin{subfigure}[b]{0.32\textwidth}\centering
		\input{./ex2_basis_functions_1.tikz}
		\caption{{LOD, node $(1/6,1/6)$}}\label{fig:ex2basisfunctions-b}
	\end{subfigure}\hfill
	\begin{subfigure}[b]{0.32\textwidth}\centering
		\input{./ex2_basis_functions_2.tikz}
		\caption{{LOD, node $(1/2,1/6)$}}\label{fig:ex2basisfunctions-c}
	\end{subfigure}
	\par\smallskip
	\begin{subfigure}[b]{0.32\textwidth}\centering
		\input{./ex2_basis_functions_3.tikz}
		\caption{{DR-LOD, node $(1/2,1/2)$}}\label{fig:ex2basisfunctions-d}
	\end{subfigure}\hfill
	\begin{subfigure}[b]{0.32\textwidth}\centering
		\input{./ex2_basis_functions_4.tikz}
		\caption{{DR-LOD, node $(1/6,1/6)$}}\label{fig:ex2basisfunctions-e}
	\end{subfigure}\hfill
	\begin{subfigure}[b]{0.32\textwidth}\centering
		\input{./ex2_basis_functions_5.tikz}
		\caption{{DR-LOD, node $(1/2,1/6)$}}\label{fig:ex2basisfunctions-f}
	\end{subfigure}
	\par\smallskip
	\begin{subfigure}[b]{0.32\textwidth}\centering
		\input{./ex2_basis_functions_6.tikz}
		\caption{{difference, node $(1/2,1/2)$}}\label{fig:ex2basisfunctions-g}
	\end{subfigure}\hfill
	\begin{subfigure}[b]{0.32\textwidth}\centering
		\input{./ex2_basis_functions_7.tikz}
		\caption{{difference, node $(1/6,1/6)$}}\label{fig:ex2basisfunctions-h}
	\end{subfigure}\hfill
	\begin{subfigure}[b]{0.32\textwidth}\centering
		\input{./ex2_basis_functions_8.tikz}
		\caption{{difference, node $(1/2,1/6)$}}\label{fig:ex2basisfunctions-i}
	\end{subfigure}
	\caption{{LOD basis functions (a)--(c), DR-LOD basis functions (d)--(f), and their difference (g)--(i) in Example~2 for the nodes $(1/2,1/2)$, $(1/6,1/6)$, and $(1/2,1/6)$ and the test coefficient at time $t=5\cdot10^3$, whose region states have not been seen during training. Each panel shows a window around the respective node. The mean relative error of the DR-LOD basis functions over the test set is $3.1\,\%$.}}\label{fig:ex2basis}
\end{figure}

\begin{figure}[htbp]
	\centering
	\begin{subfigure}[b]{0.32\textwidth}\centering
		\input{./ex2_solutions_0.tikz}
		\caption{{coefficient, $t=417$}}\label{fig:ex2solutions-a}
	\end{subfigure}\hfill
	\begin{subfigure}[b]{0.32\textwidth}\centering
		\input{./ex2_solutions_1.tikz}
		\caption{{coefficient, $t=5\cdot10^3$}}\label{fig:ex2solutions-b}
	\end{subfigure}\hfill
	\begin{subfigure}[b]{0.32\textwidth}\centering
		\input{./ex2_solutions_2.tikz}
		\caption{{coefficient, $t=10^4$}}\label{fig:ex2solutions-c}
	\end{subfigure}
	\par\smallskip
	\begin{subfigure}[b]{0.32\textwidth}\centering
		\input{./ex2_solutions_6.tikz}
		\caption{{LOD, $t=417$}}\label{fig:ex2solutions-g}
	\end{subfigure}\hfill
	\begin{subfigure}[b]{0.32\textwidth}\centering
		\input{./ex2_solutions_7.tikz}
		\caption{{LOD, $t=5\cdot10^3$}}\label{fig:ex2solutions-h}
	\end{subfigure}\hfill
	\begin{subfigure}[b]{0.32\textwidth}\centering
		\input{./ex2_solutions_8.tikz}
		\caption{{LOD, $t=10^4$}}\label{fig:ex2solutions-i}
	\end{subfigure}
	\par\smallskip
	\begin{subfigure}[b]{0.32\textwidth}\centering
		\input{./ex2_solutions_9.tikz}
		\caption{{DR-LOD, $t=417$}}\label{fig:ex2solutions-j}
	\end{subfigure}\hfill
	\begin{subfigure}[b]{0.32\textwidth}\centering
		\input{./ex2_solutions_10.tikz}
		\caption{{DR-LOD, $t=5\cdot10^3$}}\label{fig:ex2solutions-k}
	\end{subfigure}\hfill
	\begin{subfigure}[b]{0.32\textwidth}\centering
		\input{./ex2_solutions_11.tikz}
		\caption{{DR-LOD, $t=10^4$}}\label{fig:ex2solutions-l}
	\end{subfigure}
	\caption{{Heat conductivity $\overline\lambda$ (a)--(c), LOD solution (d)--(f), and DR-LOD solution (g)--(i) for the test set of Example~2 at the times $t=417$, $t=5\cdot10^3$, and $t=10^4$. The states of the $48$ AM regions, and thus the conductivity, change in time, and the DR-LOD basis is updated in each time step by evaluating the networks. The errors compared to the fine finite element solution are given in Table~\ref{tab:ex2}.}}\label{fig:ex2sol}
\end{figure}

We note that the LOD method can also deal with different boundary conditions such as Neumann- or Robin-type, which we did not consider in our experiments.

%----------------------------------------------------------
% Conclusions
%----------------------------------------------------------

\section{Conclusions and Outlook}\label{sec:conclusions}
We have presented a numerical homogenization strategy based on the LOD method using a Deep Ritz approach for computing the corrections. The goal is to leverage neural networks to overcome the bottleneck of re-computations in the case of time-dependence or uncertainties in the PDE coefficient. {The networks take the parameters of the coefficient as additional inputs, such that, once trained, they provide the corrections for any coefficient of the parametrized family by a forward pass.} {The constraint condition is enforced exactly by a fixed projection of the network output, such that the training is an unconstrained minimization of the energy. {As in the Ritz--Galerkin method, a global minimizer yields the best approximation of the LOD corrections by the networks in the energy norm, without requiring the corrections themselves, and the value of the loss directly measures the error of the network corrections in this norm.}. Combined with a normalization of the networks and a fine-tuning with L-BFGS, this yields accurate corrections without any precomputed data.} We have presented two examples to illustrate the performance of the approach{: a coefficient with a growing damage front that depends strongly nonlinearly on a single parameter, and a battery cell whose coefficient depends on $48$ parameters. In both cases, the DR-LOD solutions for coefficients that have not been seen during training deviate from the LOD solutions by less than $1\,\%$. In particular, for many parameters, the training only requires a moderate number of coefficient samples, whereas precomputing the LOD basis on a grid in the parameter space is infeasible.}   
A straight-forward extension is the amplification of the Deep Ritz method using different ANN architectures. A more technical modification would be a suitable addition of the localization parameter to the training parameters. Finally, a rigorous a posteriori analysis is still open.  

\section*{{CRediT authorship contribution statement}}
{\textbf{Mehdi Elasmi:} Conceptualization, Methodology, Software, Formal analysis, Investigation, Visualization, Writing -- original draft, Writing -- review \& editing.
\textbf{Felix Krumbiegel:} Software, Visualization, Writing -- original draft, Writing -- review \& editing.
\textbf{Roland Maier:} Conceptualization, Writing -- original draft, Writing -- review \& editing.}

\section*{{Declaration of generative AI and AI-assisted technologies.}} 
{During the preparation of this work, the authors used Claude in order to improve formulation for parts of the manuscript and for parts of the code to improve the presentation. After using this tool, the authors reviewed and edited the content as needed and take full responsibility for the content of the published article.}

\section*{Acknowledgments}

R.~Maier acknowledges funding by the Deutsche Forschungsgemeinschaft (DFG, German Research
Foundation) -- Project-ID 545165789.

{M.~Elasmi acknowledges financial support by the DFG through the Research Training Group 2218 SiMET – Simulation of Mechano-Electro-Thermal processes in Lithium-ion Batteries, project number 281041241.}

\appendix
\section{{Proofs of the results in Section~\ref{subsec:exact}}}\label{app:proofs}
\renewcommand{\thesubsection}{\Alph{section}.\arabic{subsection}}

\subsection{{Proof of Lemma~\ref{lem:projection}}}\label{app:proofProjection}
{(a) Since $\mathbf{Q}^\top\mathbf{Q}=\Id$, we have $(\mathbf{Q}\mathbf{Q}^\top)^2=\mathbf{Q}\mathbf{Q}^\top$, and the assertion follows.

(b) Recall that $\mathcal{A}^\ell=\ker\mathbf{B}=(\operatorname{range}\mathbf{B}^\top)^\perp=(\operatorname{range}\mathbf{Q})^\perp$. For $\boldsymbol{\beta}\in\R^{n_0}$, we have $\mathbf{Q}^\top P^\ell\boldsymbol{\beta}=\mathbf{Q}^\top\boldsymbol{\beta}-\mathbf{Q}^\top\boldsymbol{\beta}=\mathbf{0}$, hence $P^\ell\boldsymbol{\beta}\in\mathcal{A}^\ell$. If $\boldsymbol{\beta}\in\mathcal{A}^\ell$, then $\mathbf{Q}^\top\boldsymbol{\beta}=\mathbf{0}$ and thus $P^\ell\boldsymbol{\beta}=\boldsymbol{\beta}$.

(c) Let $\boldsymbol{\beta}\in\R^{n_0}$ and $\boldsymbol{\gamma}\in\mathcal{A}^\ell$. Then $\boldsymbol{\beta}-P^\ell\boldsymbol{\beta}=\mathbf{Q}\mathbf{Q}^\top\boldsymbol{\beta}\in\operatorname{range}\mathbf{Q}$, while $P^\ell\boldsymbol{\beta}-\boldsymbol{\gamma}\in\mathcal{A}^\ell=(\operatorname{range}\mathbf{Q})^\perp$ by~(b). Hence, the two vectors are orthogonal, and $|\boldsymbol{\beta}-\boldsymbol{\gamma}|^2=|\boldsymbol{\beta}-P^\ell\boldsymbol{\beta}|^2+|P^\ell\boldsymbol{\beta}-\boldsymbol{\gamma}|^2\ge|\boldsymbol{\beta}-P^\ell\boldsymbol{\beta}|^2$.

(d) The matrix $\mathbf{Q}\mathbf{Q}^\top$ is the orthogonal projection onto $\operatorname{range}\mathbf{B}^\top$, which is unique and coincides with $\mathbf{B}^+\mathbf{B}$.}

\subsection{{Proof of Proposition~\ref{prop:energyIdentity}}}\label{app:proofEnergy}
{Since $\mathcal{A}^\ell$ is a linear subspace, the first-order optimality condition of~\eqref{eq:discreteCorrector} reads $(\mathbf{S}\boldsymbol{\alpha}^\star-\mathbf{r})^\top\boldsymbol{\delta}=0$ for all $\boldsymbol{\delta}\in\mathcal{A}^\ell$. As $J$ is quadratic, we have for all $\boldsymbol{\delta}\in\R^{n_0}$ that
\begin{equation*}
	J(\boldsymbol{\alpha}^\star+\boldsymbol{\delta})=J(\boldsymbol{\alpha}^\star)+(\mathbf{S}\boldsymbol{\alpha}^\star-\mathbf{r})^\top\boldsymbol{\delta}+\tfrac12\,\boldsymbol{\delta}^\top\mathbf{S}\boldsymbol{\delta} .
\end{equation*}
For $\boldsymbol{\delta}:=P^\ell\boldsymbol{\beta}-\boldsymbol{\alpha}^\star$, which belongs to $\mathcal{A}^\ell$ by Lemma~\ref{lem:projection}(b), the second term vanishes, which proves~\eqref{eq:energyIdentity}. Its right-hand side is non-negative and vanishes if and only if $P^\ell\boldsymbol{\beta}=\boldsymbol{\alpha}^\star$, since $\|\cdot\|_{\mathbf{S}}$ is a norm on $\R^{n_0}$. The minimum is attained at $\boldsymbol{\beta}=\boldsymbol{\alpha}^\star$, because $P^\ell\boldsymbol{\alpha}^\star=\boldsymbol{\alpha}^\star$.}

\subsection{{Proof of Corollary~\ref{cor:lossError}}}\label{app:proofLossError}{The optimality condition with $\boldsymbol{\delta}=\boldsymbol{\alpha}^\star\in\mathcal{A}^\ell$ gives $\mathbf{r}^\top\boldsymbol{\alpha}^\star=\|\boldsymbol{\alpha}^\star\|_{\mathbf{S}}^2$, hence $J(\boldsymbol{\alpha}^\star)=-\tfrac12\|\boldsymbol{\alpha}^\star\|_{\mathbf{S}}^2$. Dividing~\eqref{eq:energyIdentity} by $|J(\boldsymbol{\alpha}^\star)|$ yields the assertion.}

\section{{Network size study}}\label{app:size}
\setcounter{figure}{0}
{In this appendix, we investigate for Example~1 how the accuracy of the networks depends on their size, i.e., the trade-off between the number of trainable parameters and the accuracy. To reduce the cost of the study, we use the coarser mesh $H=1/6$ with $\ell=1$; the coefficient, the fine mesh $h=1/36$, the training parameters, and the training procedure are the same as in Section~\ref{subsec:ex1}. We train networks of widths $16$, $32$, $64$, and $128$ and depths $2$, $4$, and $8$ with the same procedure for four representative elements (a corner element, an element at the boundary, and two interior elements), and compute the relative errors of the network corrections compared to the discrete LOD corrections~\eqref{eq:discreteCorrector} in the energy norm, cf.~Corollary~\ref{cor:lossError}, for all test parameters. The results are shown in Figure~\ref{fig:ex1size}. The error decreases with the number of weights~$N_\Theta$ at an observed rate of about $N_\Theta^{-1/2}$ when the width is increased for a fixed depth (least-squares fits of the mean error in the energy norm give the rates $0.57$, $0.64$, and $0.48$ for the depths $2$, $4$, and $8$), i.e., about first order in the width, since $N_\Theta$ grows quadratically with the width. Increasing the depth for a fixed width is more efficient, with observed rates between $0.8$ and $1.0$. For the largest networks, the rate decreases, which indicates that the error of the optimization starts to dominate over the approximation error of the networks. For instance, the mean error is $35\,\%$ for width~$16$ and depth~$2$ ($404$ weights), $1.9\,\%$ for width~$64$ and depth~$4$ ($13\,000$ weights), and $1.0\,\%$ (maximum $3.1\,\%$) for width~$128$ and depth~$8$ ($117\,000$ weights), which is the network size used in Section~\ref{subsec:ex1}. At the same time, the training time per network increases from about $3$ to $26$ minutes. For comparison, an interior element patch contains $17^2=289$ interior fine nodes, i.e., the four discrete corrections of one element for a single parameter have $1156$ degrees of freedom, and for the $40$ training parameters about $46\,000$. Note that the errors of the corrections in the energy norm are larger than the errors of the basis functions reported in Section~\ref{subsec:ex1} (for $H=1/12$, $\ell=2$), since the latter include the classical basis functions $\Lambda_j$ and are measured in a discrete $L^2$-norm.}

\begin{figure}[htbp]
	\centering
	\begin{subfigure}[b]{0.48\textwidth}\centering
		\input{./network_size_0.tikz}
		\caption{energy norm}\label{fig:ex1networksize-a}
	\end{subfigure}\hfill
	\begin{subfigure}[b]{0.48\textwidth}\centering
		\input{./network_size_1.tikz}
		\caption{$L^2$-norm}\label{fig:ex1networksize-b}
	\end{subfigure}
	\par\smallskip
	\begin{subfigure}[b]{0.48\textwidth}\centering
		\input{./network_size_2.tikz}
		\caption{energy norm against training time}\label{fig:ex1networksize-c}
	\end{subfigure}
	\caption{{Relative errors of the network corrections compared to the discrete LOD corrections in Example~1 with $H=1/6$ and $\ell=1$ for networks of different widths ($16$, $32$, $64$, $128$, indicated at the markers) and depths, for four representative elements and all test parameters: mean (solid) and maximum (dashed) error in the energy norm (a) and in the $L^2$-norm (b) against the number of weights per network, and mean error in the energy norm against the training time per network (c). The triangles indicate the rate $N_\Theta^{-1/2}$ in the number of weights~$N_\Theta$.}}\label{fig:ex1size}
\end{figure}

\bibliography{references} 
\bibliographystyle{alpha}

\end{document}